\documentclass[12pt]{article}
\usepackage{amsfonts,amsmath,amsxtra}
\usepackage{amssymb,amscd}
\usepackage[mathscr]{eucal}
\def\hybrid{\topmargin 0pt      \oddsidemargin 0pt
        \headheight 0pt \headsep 0pt
        \textwidth 16.5cm
        \textheight 23cm
        \marginparwidth 0.0in
        \parskip 5pt plus 1pt   \jot = 1.5ex}
\catcode`\@=11
\def\marginnote#1{}
\newcount\hour
\newcount\minute
\newtoks\amorpm
\hour=\time\divide\hour by60 \minute=\time{\multiply\hour by60
\global\advance\minute by-\hour}
\edef\standardtime{{\ifnum\hour<12 \global\amorpm={am}%
        \else\global\amorpm={pm}\advance\hour by-12 \fi
        \ifnum\hour=0 \hour=12 \fi
      \number\hour:\ifnum\minute<10 0\fi\number\minute\the\amorpm}}
\edef\militarytime{\number\hour:\ifnum\minute<10 0\fi\number\minute}

\def\draftlabel#1{{\@bsphack\if@filesw {\let\thepage\relax
   \xdef\@gtempa{\write\@auxout{\string
      \newlabel{#1}{{\@currentlabel}{\thepage}}}}}\@gtempa
   \if@nobreak \ifvmode\nobreak\fi\fi\fi\@esphack}
        \gdef\@eqnlabel{#1}}
\def\@eqnlabel{}
\def\@vacuum{}
\def\draftmarginnote#1{\marginpar{\raggedright\scriptsize\tt#1}}

\def\draft{\oddsidemargin -0.1truein
        \def\@oddfoot{\sl preliminary draft \hfil
        \rm\thepage\hfil\sl\today\quad\militarytime}
        \let\@evenfoot\@oddfoot \overfullrule 3pt
        \let\label=\draftlabel
        \let\marginnote=\draftmarginnote
\def\@eqnnum{{\rm (\theequation)}
\rlap{\kern\marginparsep\tt\@eqnlabel}%
\global\let\@eqnlabel\@vacuum}  }

\newcommand{\RR}{{\mathbb{R}}}
\newcommand{\CC}{{\mathbb{C}}}

\newcommand{\ZZ}{{\mathbb{Z}}}

\newfont{\Bbbb}{msbm7 scaled 1\@ptsize00}
\newcommand{\zs}{\raise-1pt\hbox{$\mbox{\Bbbb Z}$}}

scaled 1\@ptsize00

\font\sevenmsa=msam6 
\newfam\msafam
\textfont\msafam=\sevenmsa
\def\hexnumber@#1{\ifnum#1<10 \number#1\else
\ifnum#1=10 A\else\ifnum#1=11 B\else\ifnum#1=12 C\else \ifnum#1=13
D\else\ifnum#1=14 E\else\ifnum#1=15 F\fi\fi\fi\fi\fi\fi\fi}
\def\msa@{\hexnumber@\msafam}
\def\llcorner{\delimiter"4\msa@78\msa@78 }
\def\lrcorner{\delimiter"5\msa@79\msa@79 }
\mathchardef\blacktriangleright="3\msa@49
\mathchardef\blacktriangleleft="3\msa@4A \font\tenmsb=msbm10 scaled
1\@ptsize00
\newfam\msbfam
\textfont\msbfam=\tenmsb \scriptfont\msbfam=\tenmsb

\newdimen\linethick  \linethick=0.4pt
\newdimen\hboxitspace    \hboxitspace=5pt
\newdimen\vboxitspace    \vboxitspace=5pt

\def\fr#1{%
\be\new \vcenter{ \hrule height\linethick
           \hbox{\vrule width\linethick
                 \kern\hboxitspace
                 \vbox{\kern\vboxitspace
                       \hbox{$\begin{array}{c}\displaystyle#1
          \end{array}$}%
                       \kern\vboxitspace}%
                 \kern\hboxitspace
                 \vrule width\linethick}%
           \hrule height\linethick}%
\ee}

\newdimen\Squaresize \Squaresize=14pt
\newdimen\Thickness \Thickness=0.5pt

\def\Square#1{\hbox{\vrule width \Thickness
   \vbox to \Squaresize{\hrule height \Thickness\vss
      \hbox to \Squaresize{\hss#1\hss}
   \vss\hrule height\Thickness}
\unskip\vrule width \Thickness} \kern-\Thickness}

\def\Vsquare#1{\vbox{\Square{$#1$}}\kern-\Thickness}

\def\numberbysection{\@addtoreset{equation}{section}
        \def\theequation{\thesection.\arabic{equation}}}
\numberbysection

\renewcommand{\theequation}{\thesection.\arabic{equation}}
\def\titlepage{\@restonecolfalse\if@twocolumn\@restonecoltrue\onecolumn
     \else \newpage \fi \thispagestyle{empty}\c@page\z@
        \def\thefootnote{\fnsymbol{footnote}} }

\def\endtitlepage{\if@restonecol\twocolumn \else  \fi
        \def\thefootnote{\arabic{footnote}}
        \setcounter{footnote}{0}}  
\relax

\hybrid
\makeatletter
\newdimen\normalarrayskip            
\newdimen\minarrayskip               
\normalarrayskip\baselineskip \minarrayskip\jot
\newif\ifold             \oldtrue            \def\new{\oldfalse}
\def\arraymode{\ifold\relax\else\displaystyle\fi}
\def\eqnumphantom{\phantom{(\theequation)}} 
\def\@arrayskip{\ifold\baselineskip\z@\lineskip\z@
     \else
     \baselineskip\minarrayskip\lineskip1\baselineskip\fi}

\def\@arrayclassz{\ifcase \@lastchclass \@acolampacol \or
\@ampacol \or \or \or \@addamp \or
   \@acolampacol \or \@firstampfalse \@acol \fi
\edef\@preamble{\@preamble
  \ifcase \@chnum
     \hfil$\relax\arraymode\@sharp$\hfil
     \or $\relax\arraymode\@sharp$\hfil
     \or \hfil$\relax\arraymode\@sharp$\fi}}

\def\@array[#1]#2{\setbox\@arstrutbox=\hbox{\vrule
     height\arraystretch \ht\strutbox
     depth\arraystretch \dp\strutbox
width\z@}\@mkpream{#2}\edef\@preamble{\halign \noexpand\@halignto
\bgroup \tabskip\z@ \@arstrut \@preamble \tabskip\z@ \cr}%
\let\@startpbox\@@startpbox \let\@endpbox\@@endpbox
  \if #1t\vtop \else \if#1b\vbox \else \vcenter \fi\fi
  \bgroup \let\par\relax
  \let\@sharp##\let\protect\relax
  \@arrayskip\@preamble}
\def\eqnarray{\stepcounter{equation}%
              \let\@currentlabel=\theequation
              \global\@eqnswtrue
              \global\@eqcnt\z@
              \tabskip\@centering              
              \let\\=\@eqncr
              $$%
            \halign to \displaywidth  \bgroup
             \eqnumphantom \@eqnsel
      \hskip\@centering                               
    $\displaystyle  \tabskip\z@ {##}$%
    &\global\@eqcnt\@ne \hskip 2\arraycolsep
         $ \displaystyle  \arraymode{##}$\hfil
    &\global\@eqcnt\tw@ \hskip 2\arraycolsep
         $\displaystyle\tabskip\z@{##}$\hfil
         \tabskip\@centering
    &{##}\tabskip\z@\cr}
\makeatother
\newtheorem{te}{Theorem}[section]
\newtheorem{de}{Definition}[section]
\newtheorem{prop}{Proposition}[section]           
\newtheorem{cor}{Corollary}[section]
\newtheorem{lem}{Lemma}[section]

\newcommand{\beq}[1]{\begin{equation}\label{#1}}
\newcommand\eeq{\end{equation}}
\newcommand\bqa{\begin{eqnarray}}
\newcommand\eqa{\end{eqnarray}}
\def\be{\begin{eqnarray}\new\begin{array}{cc}}
\def\ee{\end{array}\end{eqnarray}}

\def\beq{\begin{equation}}
\def\eeq{\end{equation}}
\def\bse{\begin{subequations}}                
\def\ese{\end{subequations}}
\def\bp{\begin{pmatrix}}
\def\ep{\end{pmatrix}}

\def\i{\imath}

\def\stack#1#2{\raise0.7pt\hbox{$\mathrel{\mathop{#2}\limits^{#1}}$}}
\def\tr{\triangleright}
\def\tl{\triangleleft}
\def\sem{\mathsurround=0pt \raise1pt
\hbox{$\scriptscriptstyle>\!\!$}\:\!\!\tl}
\def\mes{\mathsurround=0pt \tr\!\:\!\raise0.8pt
\hbox{$\scriptscriptstyle\!\!<$}\,}
\def\]{\mathsurround=0pt ]\raise-2pt\hbox{$_\ast$}}

\def\la{\lambda}
\def\l{\lambda}

\def\<{\langle}
\def\>{\rangle}

\def\CQ{{\cal Q}}
\def\frak{\mathfrak}

\def\CO{{\cal O}}
\def\CU{{\cal U}}
\def\CZ{{\cal Z}}

\def\CH{\mathcal{H}}

\def\we{\raise-1pt\hbox{$\,\stackrel{\wedge}{,}\,$}}

\def\Tr{{\rm Tr}\,}
\def\pr {\partial}

0
\begin{document}

\footnotesize
\normalsize

\newpage

\thispagestyle{empty}

\begin{center}

\phantom.
{\hfill{\normalsize hep-th/yymmnnn}\\
\hfill{\normalsize ITEP-TH-26/07}\\
\hfill{\normalsize HMI-07-07}\\
\hfill{\normalsize TCD-MATH-07-14}\\
[10mm]\Large\bf
Baxter operator and  Archimedean Hecke algebra} \vspace{0.5cm}

\bigskip
{\large A. Gerasimov}
\\ \bigskip
{\it Institute for Theoretical and
Experimental Physics, 117259, Moscow,  Russia,} \\
{\it  School of Mathematics, Trinity
College, Dublin 2, Ireland, } \\
{\it  Hamilton
Mathematics Institute, TCD, Dublin 2, Ireland},\\
\bigskip
{\large D. Lebedev\footnote{E-mail: lebedev@itep.ru}}
\\ \bigskip
{\it Institute for Theoretical and Experimental Physics, 117259,
Moscow, Russia},\\{\it Max-Planck-Institut f\"ur Mathematik,
Vivatsgasse 7, D-53111
Bonn, Germany},\\
\bigskip
{\large S. Oblezin} \footnote{E-mail: Sergey.Oblezin@itep.ru}\\
\bigskip {\it Institute for Theoretical and Experimental Physics,
117259, Moscow, Russia},\\{\it Max-Planck-Institut f\"ur Mathematik,
Vivatsgasse 7, D-53111
Bonn, Germany}.\\
\end{center}

\vspace{0.5cm}

\begin{abstract}
\noindent

In this paper  we introduce Baxter integral $\mathcal{Q}$-operators
for finite-dimensional Lie algebras $\mathfrak{gl}_{\ell+1}$ and
$\mathfrak{so}_{2\ell+1}$. Whittaker functions corresponding to
these algebras are eigenfunctions of the $\mathcal{Q}$-operators with the eigenvalues
expressed in terms of Gamma-functions. The appearance of
the Gamma-functions is one of the manifestations of an interesting connection between
Mellin-Barnes and Givental integral representations
of Whittaker functions, which are in a sense  dual to
each other. We define a dual Baxter operator
and  derive  a family of mixed Mellin-Barnes-Givental
integral representations. Givental and Mellin-Barnes
integral representations are used to provide a short proof
of  the Friedberg-Bump and Bump conjectures for $G=GL(\ell+1)$
proved earlier by Stade.  We also identify eigenvalues of
the Baxter $\mathcal{Q}$-operator
acting on Whittaker functions with local Archimedean $L$-factors.
The Baxter $\mathcal{Q}$-operator
introduced in this paper is then
described as a particular realization of the explicitly
defined universal Baxter operator in the spherical
 Hecke algebra $\CH(G(\mathbb{R}), K)$,
$K$ being a maximal compact~subgroup~of~$G$.  Finally we
stress an  analogy between $\mathcal{Q}$-operators and certain
elements of the non-Archimedean Hecke algebra $\CH(G(\mathbb{Q}_p),G(\mathbb{Z}_p))$.

\end{abstract}

\vspace{1cm}

\clearpage \newpage



\section{Introduction}

The notion of the $\mathcal{Q}$-operator  was introduced by
Baxter  as an important tool to solve quantum integrable systems  \cite{Ba}.
These operators were initially constructed for
a particular class of  quantum integrable
systems associated with affine Lie algebras
$\widehat{\mathfrak{gl}}_{\ell+1}$ and its quantum/elliptic
generalizations. A new class of integral  $\mathcal{Q}$-operators
corresponding to  $\widehat{\mathfrak{gl}}_{\ell+1}$-Toda chain was
later proposed  by Pasquier and  Gaudin \cite{PG}.
Its generalization to  Toda chains for
other classical affine Lie algebras was proposed recently in
\cite{GLO1}, \cite{GLO2}, \cite{GLO3}.

In this paper  we introduce   integral Baxter $\mathcal{Q}$-operators
for Toda chains corresponding to the  finite-dimensional classical Lie
algebras $\mathfrak{gl}_{\ell+1}$ and $\mathfrak{so}_{2\ell+1}$.
These integral operators are closely related with the recursion
operators  in the Givental integral representation of Whittaker functions
(see \cite{Gi}, \cite{JK}  for $\mathfrak{gl}_{\ell+1}$ and \cite{GLO3} for other
classical Lie algebras).  It is well known that
$\mathfrak{g}$-Whittaker functions are common eigenfunctions
of the complete set of  mutually commuting $\mathfrak{g}$-Toda chain
quantum Hamiltonians. The quantum  Hamiltonians  arise as  projections of the
generators of the center $\CZ(\mathfrak{g})$ of the universal enveloping
algebra $\mathcal{U}(\mathfrak{g})$. One of the characteristic properties of the
introduced Baxter integral operators for a finite-dimensional classical Lie
algebra $\mathfrak{g}$  is that the corresponding
$\mathfrak{g}$-Whittaker functions are
their  eigenfunctions. Moreover, integral $\mathcal{Q}$-operators provide
 a complete set of  integral equations  defining $\mathfrak{g}$-Whittaker
functions. Similarly to the relation of the Hamiltonians with the
generators of the center $\CZ$, we construct universal Baxter
operators in a spherical Hecke algebra whose projection gives
Baxter operator for Toda chains. Other projections provide Baxter
operators for other quantum integrable systems (e.g. Sutherland
models).

The eigenvalues of the Baxter operators acting on Whittaker functions
are  expressed in terms of a product of Gamma-functions. The appearance of
the Gamma-functions  implies a close connection between
Givental and Mellin-Barnes integral representations \cite{KL1}
for  $\mathfrak{gl}_{\ell+1}$-Whittaker functions.
We discuss this relation in some detail.
Note that the representation theory interpretation \cite{GKL}
of the Mellin-Barnes integral representation uses
the Gelfand-Zetlin construction of the maximal commutative subalgebra
of $\mathcal{U}(\mathfrak{gl}_{\ell+1})$. One can guess a connection between
Mellin-Barnes and Givental representations
on a general ground by noticing that Givental diagrams
for classical Lie algebras \cite{GLO3}
are identical to Gelfand-Zetlin patterns  \cite{BZ}.
Moreover both constructions are most natural for classical Lie algebras.
In this note we discuss  a duality relation  between recursive
structures of  Givental and Mellin-Barnes integral
representations. We construct a dual version of the Baxter
$\CQ$-operator and derive a set of  relations between recursive/Baxter operators
and their duals.  We also propose  a family of
mixed Mellin-Barnes-Givental integral representations
interpolating between Mellin-Barnes and Givental integral
representations of Whittaker functions.

We  use the  Mellin-Barnes integral representation to give
simple proofs  of  Bump-Friedberg and Bump conjectures
on Archimedean factors arising in the application of the Rankin-Selberg method
to analytic continuations  of $GL(\ell+1)\times GL(\ell+1)$
and $GL(\ell+1)\times GL(\ell)$ automorphic $L$-functions.
We also  discuss a relation  with  the proofs
given by Stade \cite{St1},\,\cite{St2}.
The proof in \cite{St1},\,\cite{St2}  is based on a
recursive construction of $\mathfrak{gl}_{\ell+1}$-Whittaker
functions generalizing  the construction due to  Vinogradov and
Takhtadzhyn \cite{VT}. As it was noticed in \cite{GKLO} and is explicitly
demonstrated below, the Stade recursion basically coincides with
the   Givental recursion (see also recent  detailed discussion in
\cite{St3}). We also show  that the Bump-Friedberg and Bump
conjectures are simple consequences of the Mellin-Barnes integral
representation of $\mathfrak{gl}_{\ell+1}$-Whittaker function.

Rankin-Selberg method is a powerful tool of studying  analytic
properties of automorphic $L$-functions.
The application of the Baxter $\mathcal{Q}$-operators and closely related recursive
operators to a derivation of analytic properties of $L$-functions
using Rankin-Selberg method  is not accidental.
We remark that the eigenvalues of the  $\mathcal{Q}$-operators
acting on $\mathfrak{g}$-Whittaker functions are given by Archimedean local
$L$-factors and the  integral $\mathcal{Q}$-operators
should be naturally  considered as elements of the Archimedean Hecke
 algebra  $\CH(G(\mathbb{R}), K)$, $K$ being a maximal compact~subgroup~of~$G$.
We construct the corresponding universal Baxter operator
as an element  of the spherical Hecke algebra $\CH(G(\mathbb{R}),
K)$.
 We also describe non-Archimedean counterparts of
the universal Baxter operators as elements
of non-Archimedean Hecke algebras
$\CH(GL(\ell+1,\mathbb{Q}_p),GL(\ell+1,\mathbb{Z}_p))$.
The consideration of Archimedean and non-Archimedean universal $\CQ$-operators
on an equal footing provides  a uniform description of the
automorphic forms as their common eigenfunctions (replacing  more
traditional approach based on  the algebra of the
invariant differential operators as a substitute of  $\CH(G(\mathbb{R}), K)$).

Let us note that the connection of the Baxter operators with
Archimedean $L$-factors implies in particular that there is a hidden parameter
in the $\mathcal{Q}$-operator corresponding to a
choice of a  finite-dimensional representation of
the Langlands dual Lie algebra. In this sense,  $\mathcal{Q}$-operators considered
in this paper correspond to  standard representations of the classical
Lie algebras.  We are going to consider the
$\mathcal{Q}$-operators corresponding to more general representations
in a separate publication.

One should stress that there are various Hecke algebras relevant to
the study of the quantum Toda chains. For example for $B\subset G$
being a Borel subgroup, the Hecke algebra $\CH(G,B)$ of
$B$-biinvariant functions is closely related to  the scattering data
of quantum Toda chains \cite{STS}. The Hecke algebra $\CH(G,N)$, $N$
being the unipotent radical  of $B$ also  deserves a consideration.
Note that the representations of $\CH(G,N)$ contain  certain
information about the scattering data of  the theory and its center
is isomorphic to $\CH(G,K)$. 

Finally let us remark  that the
constructions of affine integral $\mathcal{Q}$-operators and their
eigenvalues for the action on Whittaker functions \cite{PG} 
together with the considerations of this paper  imply an
intriguing  possibility to interpret the eigenvalues of affine
$\mathcal{Q}$-operators as a kind of local  Archimedean $L$-factors.
It is natural to expect that these $L$-factors should be connected
with 2-dimensional local fields in the sense of Parshin \cite{Pa}.
We are going to discuss this fascinating possibility elsewhere.

The plan of this paper  is as follows. In Section 2
we recall Givental integral representation
for $\mathfrak{gl}_{\ell+1}$ and introduce Baxter
$\mathcal{Q}$-operator for $\mathfrak{gl}_{\ell+1}$.
In Section 3 we consider relation between Givental and Mellin-Barnes
integral representations of $\mathfrak{gl}_{\ell+1}$-Whittaker
functions and introduce dual Baxter operator.
In Section 4 we use  Mellin-Barnes integral representation
to prove Bump-Friedberg and Bump conjecture and discuss the  relation
with \cite{St1},\,\cite{St2}.
In  Section 5 we identify eigenvalues of the Baxter $\mathcal{Q}$-operator
with local Archimedean $L$-factors and construct universal Baxter
operators as elements of the spherical Hecke algebra
$\CH(G(\mathbb{R}), K)$. The main result of this  paper is given in Theorem
\ref{MainTh}. We also discuss an  analogy between
$\mathcal{Q}$-operators and  certain  elements  of the non-Archimedean
Hecke algebra $\CH(GL(\ell+1,\mathbb{Q}_p),GL(\ell+1,\mathbb{Z}_p))$.
Finally in Section 6 a generalization to $\mathfrak{so}(2\ell+1)$
is given.

{\em Acknowledgments}:  The research of AG
was  partly supported by  Science Foundation Ireland
grant and the research of SO was partially supported
by  RF President Grant MK-134.2007.1.

\section{Baxter operator for $\mathfrak{gl}_{\ell+1}$}

\subsection{Whittaker functions as matrix elements}

Let us recall two constructions of $\mathfrak{g}$-Whittaker
functions as matrix elements of infinite-dimensional representations
of $\CU(\mathfrak{g})$ and a relation of $\mathfrak{g}$-Whittaker
functions with eigenfunctions of $\mathfrak{g}$-Toda quantum chains.

Let us first describe the construction based on the Gauss decomposition.
According to Kostant \cite{Ko1}, \cite{Ko2},
 $\mathfrak{gl}_{\ell+1}$-Whittaker function can be defined as a
certain  matrix element in a principal series representation
of $G=GL(\ell+1,\mathbb{R})$. Let $\mathcal{U}(\frak{g})$ be
a universal enveloping algebra of
$\mathfrak{g}=\mathfrak{gl}_{\ell+1}$ and  $V$, $V'$  be
$\mathcal{U}(\frak{g})$-modules, dual  with respect to
a  non-degenerate invariant pairing $\<.\,,.\>: V'\times
V\to\CC$,  $\<v',Xv\>=-\<X v',v\>$ for all $v\in V$, $v'\in
V'$ and $X\in\frak{g}$.
Let $B_{-}=N_{-}AM$ and $B_+=AMN_+$ be  Langlands decompositions
of opposite Borel subgroups. Here $N_{\pm}$ are  unipotent radicals of $B_{\pm}$,
$A$ is the  identity component of the vector Cartan subgroup  and
$M$ is the intersection of the centralizer of the vector Cartan
subalgebra with the maximal compact subgroup $K\subset G$.
We will  assume that the actions of the Borel subalgebras
$\mathfrak{b}_+={\rm Lie}(B_+)$ on $V$ and
$\mathfrak{b}_-={\rm Lie}(B_-)$ on $V'$
are   integrated to the actions of the
corresponding subgroups.
Let $\chi_{\pm}\!:\mathfrak{n}_{\pm}\!\rightarrow\CC$ be the characters
of $\mathfrak{n}_{\pm}$ defined by $\chi_+(e_{i}):=-1$ and
$\chi_-(f_{i}):=-1$ for all $i=1,\ldots,\ell$.  A vector $\psi_R\in
V$ is called a Whittaker vector with respect to
 $\chi_+$ if
\be \hspace{2cm} e_i\psi_R=-\psi_{R}\,,
 \hspace{1cm}(i=1,\ldots,\ell), \ee and a vector
$\psi_L\in V'$ is called a Whittaker vector with respect to $\chi_-$
if \be \hspace{2cm} f_{i} \psi_L=-\psi_L\,, \hspace{1cm}
(i=1,\ldots,\ell). \ee One defines a  Whittaker model $\mathcal{V}$
as a space of functions on $G$ such that $f(ng)=\chi_{N_+}(n)f(g)$,
$n\in N_+$, $\chi_{N_+}(n)=\chi_+(\log n)$.
$\mathcal{U}(\mathfrak{g})$-module admits a Whittaker model with
respect to the character $\chi$ if it is equivalent to a
sub-representation of $\mathcal{V}$.

Let $\mathcal{V}_{\underline{\lambda}}={\rm
Ind}_B^G\chi_{\underline{\lambda}}$ be a principal series
representation  of $G$ induced from the generic character
$\chi_{\underline{\lambda}}$ of $B$ trivial on $N\subset B$.
It is realized in the space of
functions $f\in C^{\infty}(G)$ satisfying equation
$$f(bg)=\chi_{\underline{\lambda}}(b)f(g),$$  where $b\in B$.
The action of $G$ is given by the
right action $\pi_{\lambda}(g)f(x)=f(xg^{-1}).$ We will be
interested in the infinitesimal form $Ind_{U(\frak b)}^{U(\frak
g)}\chi_{\underline{\lambda}}$ of this representation given by
$$(Xf)(g)=\frac{d}{dt}f(ge^{-tX})|_{t\rightarrow \infty}.$$

Define $(\mathfrak{g},B)$-module as $\mathfrak{g}$-module such that
the action of the Borel subalgebra $\mathfrak{b}\subset
\mathfrak{g}$  is integrated to the action of the Borel subgroup
$B$, $\mathfrak{b}={\rm Lie}(B)$. Consider an  irreducible
$(\mathfrak{g},B)$-submodule
$\mathcal{V}^{(0)}_{\underline{\lambda}}$ of
$\mathcal{V}_{\underline{\lambda}}$ given  by the Schwartz space
$\mathcal{S}(N_-)$ of functions on $N_-$ exponentially decreasing
at the infinity with all their derivatives. This
$(\mathfrak{g},B)$-module always admits a Whittaker model. Below we
will denote by $\psi_L$, $\psi_R$ the Whittaker vectors in
$\mathcal{V}^{(0)}_{\underline{\lambda}}$ and its dual. Following
Kostant \cite{Ko1},\cite{Ko2} ( see also  \cite{Et} for a recent
discussion) we define a  $\mathfrak{g}$-Whittaker function  in terms
of the invariant pairing of  Whittaker modules as follows
 \bqa\label{pairing}
 \Psi^{\mathfrak{g}}_{\underline{\lambda}} (x)=e^{-
\<\rho,x\>}\<\psi_L\,, \pi_{\underline{\lambda}}(e^{h_x})
\,\psi_R\>,\, \qquad \qquad x\in \mathfrak{h}, \eqa  where
 $h_x:=\sum\limits_{i=1}^{\ell}\<\omega_i,x\>\,h_i$,
$\omega_i$ is a bases of fundamental weights of $\mathfrak{g}$,
$\rho=1/2\sum_{\alpha>0}\alpha$ and
$\pi_{\underline{\lambda}}(e^{h_x})$ is an action of $e^{h_x}$ in
the representation $\mathcal{V}_{\underline{\lambda}}$. It was shown
in \cite{Ko1} that $\mathfrak{g}$-Whittaker function is a common
eigenfunction of a complete  set  of commuting Hamiltonians of
$\mathfrak{g}$-Toda chain. A complete set of commuting Hamiltonians
of the $\mathfrak{g}$-Toda chain is generated by the differential
operators $\CH_k\in{\rm Diff} (\mathfrak{h})$, $k=1,\cdots, \ell$ on
the Cartan subalgebra $\mathfrak{h}$
 defined in terms  of  the  generators  $\{c_k\}$ of the center
$\mathcal{Z}(\mathfrak{g})\subset \mathcal{U}(\mathfrak{g})$  as
follows: \bqa \label{hamdef} \CH_k
\Psi^{\mathfrak{g}}_{\underline{\lambda}}(x)=e^{-\<\rho,x\>}
\<\psi_L\,,\pi_{\underline{\lambda}}(e^{h_x})\,c_k\,\psi_R\>.\qquad
\eqa Let $c_2$ be a quadratic generator of
$\mathcal{Z}(\mathfrak{g})$ (Casimir element) \bqa\label{seccas}
c_2=\frac{1}{2}\sum\limits_{i,j=1}^{\ell}c_{ij}h_ih_j+\frac{1}{2}
\sum\limits_{\alpha\in
R_+}(e_{\alpha}f_{\alpha}+f_{\alpha}e_{\alpha}),\eqa where  $R_+$ is
a set of all positive roots, $\|c_{ij}\|=\|d_id_j(b^{-1})_{ij}\|$
and $\|b_{ij}\|=\|d_i a_{ij}\|$ is a symmetrisation of the Cartan
matrix $\|a_{ij}\|$.   Let $\{\epsilon_i\}$ be an orthogonal bases
$(\epsilon_i,\epsilon_j)=\delta_{ij}$ in $\mathfrak{h}$ and
$x=\sum_{i=1}^{\ell} x_i\epsilon_i$ be  a decomposition of $x\in
\mathfrak{h}$ in this bases. Then the projection  (\ref{hamdef}) of
(\ref{seccas}) gives the following quadratic Hamiltonian operator of
$\mathfrak{g}$-Toda chain (see e.g. \cite{RSTS})
\be\CH^{\mathfrak{g}}_2=-\frac{1}{2}\sum_{i=1}^{\ell}
\frac{\partial^2}{{\partial x_i}^2} +\sum_{i=1}^\ell d_i
e^{\<\alpha_i,x\>}.\ee For $\mathfrak{g}=\mathfrak{gl}_{\ell+1}$ one
has \bqa\label{pairgl}
\Psi_{\underline{\lambda}}^{\mathfrak{gl}_{\ell+1}}(\underline{x})
=e^{- \sum_{i=1}^{\ell+1} x_i\rho_i}\<\psi_L\,,
\pi_{\underline{\lambda}}(e^{\sum_{i=1}^{\ell+1} x_i E_{i,i}})
\,\psi_R\>,\eqa where $\rho_j=\frac{\ell}{2}+1-j,j=1,\ldots,\ell+1$
are the components of $\rho$ in the standard basis of
$\RR^{\ell+1}$, $\underline{x}=(x_1,\ldots,x_{{\ell}+1})$ and
$E_{i,j}$ are the standard generators of
$\mathcal{U}(\mathfrak{gl}_{\ell+1})$. The linear and quadratic
Hamiltonians in this case are given by \bqa
\CH_1^{\mathfrak{gl}_{\ell+1}}=-\imath \sum\limits_{i=1}^{\ell+1}
\frac{\partial}{\partial x_i},  \eqa \bqa
\tilde{\CH}_2^{\mathfrak{gl}_{\ell+1}}=-\frac{1}{2}\sum\limits_{i=1}^{\ell+1}
\frac{\partial^2}{{\partial x_i}^2}+ \sum\limits_{i=1}^{\ell}
e^{x_{i}-x_{i+1}}.  \eqa Let us introduce a  generating function for
$\mathfrak{gl}_{\ell+1}$-Toda  chain Hamiltonians as
\bqa\label{genfunctgl} t^{\mathfrak{gl}_{\ell+1}}(\la)=
\sum_{j=1}^{\ell+1}\,(-1)^j\la^{\ell+1-j}\CH_j^{\mathfrak{gl}_{\ell+1}}(x,\pr_x),
\eqa where
$\tilde{\CH}_2^{\mathfrak{gl}_{\ell+1}}=\frac{1}{2}({\CH}_1^{\mathfrak{gl}_{\ell+1}})^2
-\CH_2^{\mathfrak{gl}_{\ell+1}}$. Then the
$\mathfrak{gl}_{\ell+1}$-Whittaker function  satisfies the following
equation \bqa \label{eigeneq}
t^{\mathfrak{gl}_{\ell+1}}(\la)\,\,\,\Psi^{\mathfrak{gl}_{\ell+1}}_{\underline{\la}}
 (\underline{x})= \prod_{j=1}^{\ell+1}(\la-\la_j)\,\,\,
\Psi^{\mathfrak{gl}_{\ell+1}}_{\underline{\la}}
(\underline{x}), \eqa
where $\underline{\la}=(\la_1,\cdots ,\la_{\ell+1})$ and
$\underline{x}=(x_1,\ldots,x_{{\ell}+1})$.

The appropriately normalized $\mathfrak{gl}_{\ell+1}$-Whittaker function (\ref{pairing})
is a  solution of the equations (\ref{eigeneq})
invariant with respect to the actions of
the Weyl group $W=S_{\ell+1}$ given by
$s:\lambda_i\to \lambda_{s(i)}$, $s\in W$.
The $W$-invariant $\mathfrak{gl}_{\ell+1}$- Whittaker functions provide a bases of
$W$-invariant functions in $\mathbb{R}^{\ell+1}$ (see e.g.
\cite{STS}, \cite{KL2}).

\begin{te} For the properly normalized $W$-invariant
$\mathfrak{gl}_{\ell+1}$-Whittaker functions the following
orthogonality and completeness relations hold \bqa\label{xorth}
\int_{\mathbb{R}^{\ell+1}}
\overline{\Psi}^{\mathfrak{gl}_{\ell+1}}_{\underline{\la}}
(\underline{x})
\,\Psi^{\mathfrak{gl}_{\ell+1}}_{\underline{\la}'}(\underline{x})
\,\,\prod_{j=1}^{\ell+1} d{x_j}=
\frac{1}{(\ell+1)!\,\,\mu^{(\ell+1)}(\underline{\la})}\,\,\,\sum_{w\in
W}\delta^{(\ell+1)}(\underline{\la}-w(\underline{\la}')), \eqa
 \bqa\label{lorth}
\int_{\mathbb{R}^{\ell+1}}
\overline{\Psi}^{\mathfrak{gl}_{\ell+1}}_{\underline{\la}}(\underline{x})\,
\Psi^{\mathfrak{gl}_{\ell+1}}_{\underline{\la}}(\underline{y})
\mu^{(\ell+1)}(\underline{\la})\prod_{j=1}^{\ell+1} d{\la_j}=
\delta^{(\ell+1)}(\underline{x}-\underline{y}),
 \eqa where
 \bqa\label{measure} \mu^{(\ell+1)}(\underline{\la})=\frac{1
}{(2\pi)^{\ell+1} (\ell+1)!}\,\prod_{j\neq k} \frac{1}{\Gamma
(\imath \lambda_k-\imath\lambda_j)}.\eqa
\end{te}

There exists another  construction of
$\mathfrak{gl}_{\ell+1}$-Whittaker functions that uses a pairing
of the spherical vector (i.e. a vector invariant with respect to the
maximal compact subgroup $K=SO(\ell+1,\RR)$ of $GL(\ell+1,\RR)$) and
a Whittaker vector (see e.g. \cite{J}, \cite{Ha}). Consider the
following function
\bqa\label{SpherW}
\widetilde{\Psi}^{\mathfrak{gl}_{\ell+1}}_{\underline{\la}}(g)
= e^{-\rho(g)}\<\phi_K,\pi_{\underline{\la}}(g) \,\psi_R\>, \eqa
where $\rho(g)$ is  given by $\rho(kan)=\<\rho,\log a\>$,  $\phi_K$
is a spherical vector in $\mathcal{V}_{\underline{\la}}$ \bqa
\phi_K(bgk)=\chi_{\underline{\la}}(b)\phi_K(g),\qquad \qquad
k\in K,\,b\in B_+. \eqa The function
$\widetilde{\Psi}^{\mathfrak{gl}_{\ell+1}}_{\underline{\la}}(g)$ defined by
  (\ref{SpherW}) satisfies the functional equation
\bqa \widetilde{\Psi}^{\mathfrak{gl}_{\ell+1}}_{\underline{\la}}(kgn)
=\tilde{\chi}_{N_-}(n)\,\tilde{
\Psi}^{\mathfrak{gl}_{\ell+1}}_{\underline{\la}}(g),\qquad k\in K,\,n\in N_-,
\eqa where
$\tilde{\chi}_{N_-}(n)=\exp(2\sum_{j=1}^{\ell}n_{j+1,j})$. Thus
(\ref{SpherW})  descends to a function  on the space $A$ of the
diagonal matrices $a={\rm diag}(e^{\tilde{x}_1},\ldots,
e^{\tilde{x}_{\ell+1}})$ entering the Iwasawa decomposition
$KAN_-\to GL(\ell+1,\RR)$. We fix a normalization of the matrix
element so that the function (\ref{SpherW}) is $W$-invariant.  The
resulting function on $A$ is related to  the
$\mathfrak{gl}_{\ell+1}$-Whittaker function (\ref{pairing}) by a
simple redefinition of the variables.

\begin{lem}\label{Iwarep}
The following relation between
$\Psi^{\mathfrak{gl}_{\ell+1}}_{\underline{\lambda}}
(\underline{x})$ and
$\widetilde{\Psi}^{\mathfrak{gl}_{\ell+1}}_{\underline{\tilde{\lambda}}}(\underline{\tilde{x}})$
 holds
\bqa
\widetilde{\Psi}^{\mathfrak{gl}_{\ell+1}}_{\underline{\tilde{\lambda}}}
(\underline{\tilde{x}})=\Psi^{\mathfrak{gl}_{\ell+1}}_{\underline{\lambda}}
(\underline{x}), \eqa where  $\underline{\tilde{x}}
=(\tilde{x}_1,\ldots, \tilde{x}_{\ell+1})$,
$\underline{\tilde{\lambda}} =(\tilde{\lambda}_1,\ldots,
\tilde{\lambda}_{\ell+1})$ are expressed through
$\underline{x}=(x_1,\ldots, x_{\ell+1})$,
$\underline{\lambda}=(\lambda_1,\ldots, \lambda_{\ell+1})$ as
follows
$$
\tilde{x}_j=\frac{1}{2}x_j, \qquad \tilde{\lambda}_j=2\lambda_j.
$$
\end{lem}

\subsection{Recursive and Baxter operators}

The following integral representation for
$\mathfrak{gl}_{\ell+1}$-Whittaker function was introduced by
Givental \cite{Gi} (see also \cite{JK}).

\begin{te}\label{theoremone}
$\mathfrak{gl}_{\ell+1}$-Whittaker functions (\ref{pairgl}) admit an
integral representation \bqa\label{giv}
\Psi_{\lambda_1,\ldots,\lambda_{\ell+1}}^{\mathfrak{gl}_{\ell+1}}
(x_1,\ldots,x_{\ell+1})= \int_{\RR^{\frac{\ell(\ell+1)}{2}}}
\prod_{k=1}^{\ell}\prod_{i=1}^kdx_{k,i}\,\,
e^{\mathcal{F}^{\mathfrak{gl}_{\ell+1}}(x) }, \eqa where
\bqa\label{intrep}
\mathcal{F}^{{\mathfrak{gl}}_{\ell+1}}(x)=\imath\sum\limits_{k=1}^{\ell+1}
\lambda_k\Big(\sum\limits_{i=1}^{k}
x_{k,i}-\sum\limits_{i=1}^{k-1}x_{k-1,i}\Big)-\\
\nonumber \sum\limits_{k=1}^{\ell} \sum\limits_{i=1}^{k}
\Big(e^{x_{k+1,i}-x_{k,i}}+e^{x_{k,i}-x_{k+1,i+1}}\Big),\eqa
and  $x_i:=x_{\ell+1,i},\,\,\,i=1,\ldots,\ell+1$.
 \end{te}

The interpretation of the Givental integral formula as
a matrix element (\ref{pairgl}) was
first obtained in \cite{GKLO},  where it was also noted
that the integral representation (\ref{giv})  of
$\mathfrak{gl}_{\ell+1}$-Whittaker function  has a recursive
structure over the rank $\ell$ of the Lie algebra $\mathfrak{gl}_{\ell+1}$.

\begin{cor}  The following integral operators
$Q^{\mathfrak{gl}_{k+1}}_{\mathfrak{gl}_{k}}$ provide a recursive
construction of $\mathfrak{gl}_{\ell+1}$-Whittaker functions:
 \bqa\label{recursgl}
\Psi^{\mathfrak{gl}_{\ell+1}}_{\lambda_1,\ldots,\lambda_{\ell+1}}
(\underline{x}_{\ell+1})= \int_{\RR^{\ell}}
\prod_{i=1}^{\ell}dx_{\ell,i}\,\,
Q^{\mathfrak{gl}_{\ell+1}}_{\mathfrak{gl}_{\ell}}
(\underline{x}_{\ell+1},\underline{x}_{\ell}|\lambda_{\ell+1})
\Psi_{\lambda_1,\ldots,\lambda_{\ell}}^{\mathfrak{gl}_{\ell}}
(\underline{x}_{\ell}),\eqa \bqa\label{QRecA}
Q^{\mathfrak{gl}_{\ell+1}}_{\mathfrak{gl}_{\ell}}
(\underline{x}_{\ell+1},\underline{x}_{\ell}|\lambda_{\ell+1})=\eqa
\bqa
\hspace{-0.3cm}=\exp\left\{\i\lambda_{\ell+1}\Big(\sum_{i=1}^{\ell+1}
x_{\ell+1,i}- \sum_{i=1}^{\ell}x_{\ell,i}\Big)-
\sum_{i=1}^{\ell}\Big(e^{x_{\ell+1,i}-x_{\ell,i}}+e^{x_{\ell,i}-x_{\ell+1,i+1}}\Big)
\right\},\nonumber\eqa where
$\underline{x}_{k}=(x_{k,1},\ldots,x_{k,k})$ and   we assume that
$Q^{\mathfrak{gl}_{1}}_{\mathfrak{gl}_{0}}(x_{11}|\la_1)
=e^{\imath\lambda_1x_{1,1}}.$
\end{cor}

\begin{de} Baxter operator $\mathcal{Q}^{\mathfrak{gl}_{\ell+1}}(\la)$
 for $\mathfrak{gl}_{\ell+1}$ is an
integral operator with the kernel
\bqa
\mathcal{Q}^{\mathfrak{gl}_{\ell+1}}
(\underline{x},\,\underline{y}|\,\la)=
\exp\Big\{\,\imath\la\sum_{i=1}^{\ell+1}(x_i-y_i)\,-
\sum_{k=1}^{\ell}\Big(e^{x_i-y_i}+e^{y_i-x_{i+1}}\Big)\,-\,
e^{x_{\ell+1}-y_{\ell+1}} \Big\},\eqa
 where we assume $x_{i}:=x_{\ell+1,i}$ and $y_{i}:=y_{\ell+1,i}$.
\end{de}
Note that the Baxter operator defined above  is non-trivial
even for $\mathfrak{gl}_1$.

\begin{te}
Baxter operator $\mathcal{Q}^{\mathfrak{gl}_{\ell+1}}(\la)$
satisfies the following identities
\bqa\label{firstpr}
\mathcal{Q}^{\mathfrak{gl}_{\ell+1}}(\la)\cdot
 \mathcal{Q}^{\mathfrak{gl}_{\ell+1}}(\la')=
\mathcal{Q}^{\mathfrak{gl}_{\ell+1}}(\la')\cdot
\mathcal{Q}^{\mathfrak{gl}_{\ell+1}} (\la), \eqa
\bqa\label{secondonepr} \mathcal{Q}^{\mathfrak{gl}_{\ell+1}}(\gamma)
Q^{\mathfrak{gl}_{\ell+1}}_{\mathfrak{gl}_{\ell}}
(\la)=\Gamma(\imath \gamma -\imath \la)\,\,
Q^{\mathfrak{gl}_{\ell+1}}_{\mathfrak{gl}_{\ell}} (\la)
\mathcal{Q}^{\mathfrak{gl}_{\ell}}(\gamma), \eqa
\bqa\label{secondpr} \mathcal{Q}^{\mathfrak{gl}_{\ell+1}}(\la)\cdot
T^{\mathfrak{gl}_{\ell+1}}(\la')=
T^{\mathfrak{gl}_{\ell+1}}(\la')\cdot
\mathcal{Q}^{\mathfrak{gl}_{\ell+1}}(\la), \eqa \bqa\label{thirdpr}
\mathcal{Q}^{\mathfrak{gl}_{\ell+1}}(\la-\imath)=\imath^{\ell+1}
\,\,T^{\mathfrak{gl}_{\ell+1}}(\la)\,\,\mathcal{Q}^{\mathfrak{gl}_{\ell+1}}(\la),
\eqa where \bqa\label{Genar}
T^{\mathfrak{gl}_{\ell+1}}(\underline{x},\underline{y}|\la)=
t^{\mathfrak{gl}_{\ell+1}}(\underline{x},\pr_{\underline{x}}|\la)\delta^{\ell+1}(x-y),
\eqa \bqa t^{\mathfrak{gl}_{\ell+1}}
(\underline{x},\pr_{\underline{x}}|\la) =\sum_{j=1}^{\ell+1}(-1)^j
\la^{\ell+1-j}\,\CH^{\mathfrak{gl}_{\ell+1}}_j
(\underline{x},\pr_{\underline{x}}). \eqa
\end{te}
{\it Proof}. The commutativity of $\mathcal{Q}$-operators \bqa
\int_{\mathbb{R}^{\ell+1}}\,
\CQ^{\mathfrak{gl}_{\ell+1}}(\underline{y},\underline{x}|\la)\,\,
\CQ^{\mathfrak{gl}_{\ell+1}}(\underline{x},\underline{z}|\la')\,\,\prod_{j=1}^{\ell+1}d{x_j}
=\\\int_{\mathbb{R}^{\ell+1}}\,
\CQ^{\mathfrak{gl}_{\ell+1}}(\underline{y},\underline{x}|\la')\,\,
\CQ^{\mathfrak{gl}_{\ell+1}}(\underline{x},\underline{z}|\la)\,\,\prod_{j=1}^{\ell+1}d{x_j},
\eqa is proved using the following change of  variables $x_i$
\be x_1\longmapsto-x_1+z_1+\ln\Big(e^{y_1}+e^{z_2}\Big), &\\
\nonumber x_i\longmapsto-x_i-\ln\Big(e^{-y_{i-1}}+e^{-z_i}\Big)+
\ln\Big(e^{y_i}+e^{z_{i+1}}\Big), & 1<i\leq\ell,\\
\nonumber x_{\ell+1}\longmapsto-x_{\ell+1}+y_{\ell+1}-
\ln\Big(e^{-y_\ell}+e^{-z_{\ell+1}}\Big). &\ee
The proof of (\ref{secondonepr}) is similar to the proof of
the commutativity (\ref{firstpr}).
The commutation relations (\ref{thirdpr}) and the difference equation
(\ref{thirdpr}) then easily follow from (\ref{secondonepr})
and (\ref{xorth}), (\ref{lorth}) $\Box$

\begin{cor}\label{color}
 The following relation holds
\bqa\label{eigenprop}
\int_{\RR^{\ell+1}}\,\prod_{i=1}^{\ell+1}\,dx_{i}\,\,
\mathcal{Q}^{\mathfrak{gl}_{\ell+1}}(\underline{y},
\,\underline{x}|\,\gamma)\,
\Psi^{\mathfrak{gl}_{\ell+1}}_{\underline{\lambda}}(\underline{x})\,=\,
\prod_{i=1}^{\ell+1}\Gamma(\imath \gamma -\imath\lambda_i)\,\,
\Psi^{\mathfrak{gl}_{\ell+1}}_{\underline{\lambda}}(\underline{y}),\eqa
where $\underline{x}=(x_1,\ldots, x_{\ell+1})$,
$\underline{y}=(y_1,\ldots, y_{\ell+1})$
 and $\underline{\la}=(\la_1,\ldots, \la_{\ell+1})$.
\end{cor}

Finally let us provide an expression for the kernel of the Baxter
$\CQ$-operator in the parametrization naturally arising in the
construction of $\mathfrak{gl}_{\ell+1}$-Whittaker
 functions using Iwasawa decomposition (see (\ref{SpherW}) and  Lemma \ref{Iwarep}).
Let $\tilde{\CQ}^{\mathfrak{gl}_{\ell+1}}
(\tilde{\underline{x}},\tilde{\underline{y}}|\tilde{\l})$ be defined
by
$$\tilde{\CQ}^{\mathfrak{gl}_{\ell+1}}
(\tilde{\underline{x}},\tilde{\underline{y}}|\tilde{\l})=2^{\ell+1}\exp\Big\{\imath\tilde{\l}
\sum_{i=1}^{\ell+1} (\tilde{x}_i-\tilde{y}_i)-
\sum_{k=1}^{\ell}\Big(e^{2(\tilde{x}_k-\tilde{y}_k)}+e^{2(\tilde{y}_k-\tilde{x}_{k+1})}\Big)-
 e^{2(\tilde{x}_{\ell+1}-\tilde{y}_{\ell+1})}\Big\}.$$

\begin{prop}
 The following relation holds
\bqa\label{Iwaeigenprop}
\int_{\RR^{\ell+1}}\,\prod_{i=1}^{\ell+1}\,d\tilde{x}_{i}\,\,
\widetilde{\CQ}^{\mathfrak{gl}_{\ell+1}}(\tilde{\underline{y}},
\,\tilde{\underline{x}}|\,\tilde{\gamma})\,
\tilde{\Psi}^{\mathfrak{gl}_{\ell+1}}_{\underline{\tilde{\lambda}}}
(\underline{\tilde{x}})\,=\, \prod_{i=1}^{\ell+1}
\Gamma\bigl(\frac{\i\tilde{\gamma}-\i\tilde{\la_i}}{2}\bigr)\,\,
\tilde{\Psi}^{\mathfrak{gl}_{\ell+1}}_{\underline{\tilde{\lambda}}}
(\underline{\tilde{y}})\quad\eqa
\end{prop}

\section{Givental versus  Mellin-Barnes  integral representations}

An important property of the Givental integral representation
is its recursive structure  with respect to the rank of the
Lie algebra. There is another integral representation \cite{KL1} for
$\mathfrak{gl}_{\ell+1}$-Whittaker functions generalizing
Mellin-Barnes integral representation for low ranks. This
representation also has a recursive structure. Its interpretation in terms of
representation theory uses  the
Gelfand-Zetlin construction of a maximal commutative subalgebra
in $\mathcal{U}(\mathfrak{gl}_{\ell+1})$  \cite{GKL}.
In this section we compare recursive structures of Givental
and Mellin-Barnes representations and demonstrate
that these two integral representations  should be considered as dual
to each other. We propose the construction of the dual Baxter operator
based on  Mellin-Barnes  integral representations.
We also construct a family of new integral
representations interpolating between Givental and  Mellin-Barnes
representations. Finally we introduce a symmetric recursive
construction of $\mathfrak{gl_{\ell+1}}$-Whittaker functions
such that the corresponding recursive operator is expressed through the Baxter and
dual Baxter operators. The Givental and Mellin-Barnes integral
representations are then obtained from the symmetric integral
representations by simple manipulations.

Let us first recall the Mellin-Barnes integral
representation of   $\mathfrak{gl}_{\ell+1}$-Whittaker functions.

\begin{te}  The following integral representation of
 $\mathfrak{gl}_{\ell+1}$-Whittaker function holds
\be\label{wf6} \Psi^{\mathfrak{gl}_{\ell+1}}
_{\underline{\la}}(\underline{x})=\int\limits_{\cal
S}\prod_{n=1}^{\ell}
\frac{\prod\limits_{k=1}^n\prod\limits_{m=1}^{n+1}
\Gamma(\imath\gamma_{n k}-\imath \gamma_{n+1,m})} {(2\pi)^n
n!\prod\limits_{s\neq p} \Gamma(\imath\gamma_{ns}-\imath
\gamma_{np})} e^{-\imath \sum\limits_{n=1}^{\ell+1}\,
\sum\limits_{j=1}^{\ell+1} (\gamma_{nj}-\gamma_{n-1,j})x_n}
\prod_{\stackrel{\scriptstyle n=1}{j\leq n}}^{\ell}d\gamma_{nj}\,,
\ee where
$\underline{\la}=(\la_1,\ldots,\la_{\ell+1}):=(\gamma_{\ell+1,1},\ldots,
\gamma_{\ell+1,\ell+1})$, $\underline{x}=(x_1,\ldots,x_{\ell+1})$
and the domain of  integration  $\!{\cal S}\!$ is defined by the
conditions
 $\!\max_{j}\{{\rm Im}\,\gamma_{kj}\}<
\min_m\{{\rm Im}\,\gamma_{k+1,m}\}\!$ for all $ k=1,\ldots,\ell$.
Recall that we  assume   $\gamma_{nj}=0$ for $j>n$.
\end{te}

\begin{cor}\label{corone}
The following recursive relation holds:
\bqa\label{rec}
\Psi^{\mathfrak{gl}_{\ell+1}}_{\underline{\gamma}_{\ell+1}}
(x_1,\ldots, x_{\ell+1})=\int\limits_{{\cal
    S}_{\ell}}\,\,\,\widehat{Q}^{\mathfrak{gl}_{\ell+1}}_{\mathfrak{gl}_{\ell}}
(\underline{\gamma}_{\ell+1},\underline{\gamma}_{\ell}|x_{\ell+1})
\,\,\,\,\Psi^{\mathfrak{gl}_{\ell}}_{\underline{\gamma}_{\ell}}
(x_1,\ldots ,x_{\ell})\,\,\mu^{(\ell)}(\underline{\gamma}_{\ell})
\,\, \prod\limits_{ j=1}^{\ell}d\gamma_{\ell,j}\,, \eqa where
\be\label{fourc}
\widehat{Q}^{\mathfrak{gl}_{\ell+1}}_{\mathfrak{gl}_{\ell}}
(\underline{\gamma}_{\ell+1},\underline{\gamma}_{\ell}|x_{\ell+1})
=e^{-\imath(\sum\limits_{j=1}^{\ell+1}
\gamma_{\ell+1,j}-\sum\limits_{k=1}^{\ell}\gamma_{\ell,k})x_{\ell+1}}
\prod_{k=1}^{\ell}\prod_{m=1}^{\ell+1} \Gamma(\imath
\gamma_{\ell,k}-\imath \gamma_{\ell+1,m}), \ee the measure
$\mu^{(\ell)}(\underline{\gamma}_{\ell})$ is defined by
(\ref{measure}) and $\underline{\gamma}_{k}=(\gamma_{k,1},\ldots,
\gamma_{k,k})$. We imply
$\Psi^{\mathfrak{gl}_{1}}_{\gamma_{1,1}}(x_1)=e^{-\i
\gamma_{1,1}x_{1}}$. The domain  of integration ${\cal S}_{\ell}$ is
defined by the conditions  $\max_j\{{\rm Im}\gamma_{\ell,j}\}<$\\$
\min_{m}\{{\rm Im}\gamma_{\ell+1,m}\}.$
\end{cor}

We call the integral operator
Corollary \ref{corone} the Mellin-Barnes recursive operator.

Let us stress that the recursive structure of the Mellin-Barnes
integral representation of $\mathfrak{gl}_{\ell+1}$-Whittaker
functions is   dual to that of the Givental integral representation.
Indeed, Givental recursive operator
$Q^{\mathfrak{gl}_{\ell+1}}_{\mathfrak{gl}_{\ell}}$ depends on an
additional  ``spectral'' variable $\la_{\ell+1}$ and acts in the
space of functions of the ``coordinate'' variables $\underline{x}$,
while the dual Mellin-Barnes recursive operator
$\widehat{Q}^{\mathfrak{gl}_{\ell+1}}_{\mathfrak{gl}_{\ell}}$
depends on additional ``coordinate'' variable $x_{\ell+1}$ and acts
in the space of functions of the ``spectral''  variables
$\underline{\gamma}$. Using the orthogonal and completeness
relations (\ref{xorth}), (\ref{lorth}) one can show that these two
operators are related by a conjugation by the integral operator with
the kernel
$\Psi^{\mathfrak{gl}_{\ell}}_{\underline{\gamma}_{\ell}}(\underline{x}_{\ell})$.

\begin{prop} The following integral representation
for the kernel of the recursive operator
$\widehat{Q}^{\mathfrak{gl}_{\ell+1}}_{\mathfrak{gl}_{\ell}}$ holds
\bqa\label{fourfour}
\widehat{Q}^{\mathfrak{gl}_{\ell+1}}_{\mathfrak{gl}_{\ell}}
(\underline{\gamma}_{\ell+1},\underline{\gamma}_{\ell}|x_{\ell+1,\ell+1})
=\int_{\RR^{\ell}}\,\prod_{j=1}^{\ell}d{x}_{\ell+1,j}\,\,
\overline{\Psi}^{\mathfrak{gl}_{\ell}}_{\underline{\gamma}_{\ell}}
(\underline{x}'_{\ell+1})\,\,\,
\Psi^{\mathfrak{gl}_{\ell+1}}_{\underline{\gamma}_{\ell+1}}
(\underline{x}_{\ell+1})=\nonumber \eqa \bqa
=\int_{\RR^{\ell}}\,\,\prod_{j=1}^{\ell}d{x}_{\ell+1,j}
\prod_{j=k}^{\ell}d{x}_{\ell,k}\,
\overline{\Psi}^{\mathfrak{gl}_{\ell}}_{\underline{\gamma}_{\ell}}
(\underline{x'}_{\ell+1})\,\,\,
Q^{\mathfrak{gl}_{\ell+1}}_{\mathfrak{gl}_{\ell}}
(\underline{x}_{\ell+1},\underline{x}_{\ell}|\gamma_{\ell+1})
\Psi^{\mathfrak{gl}_{\ell}}_{\underline{\gamma}_{\ell}}
(\underline{x}_{\ell}), \eqa where
$\underline{x}_{k}=(x_{k,1},\ldots ,x_{k,k})$,
$\underline{x}'_{k}=(x_{k,1},\ldots ,x_{k,k-1})$.
\end{prop}

In  view of the above duality for the recursive operators
it is natural to
introduce an operator dual to the Baxter $\CQ$-operator.

\begin{de} The dual Baxter operator
  $\widehat{\CQ}^{\mathfrak{gl}_{\ell+1}}(z)$
is an integral operator with the  kernel \bqa\label{dualBax}
\widehat{\CQ}^{\mathfrak{gl}_{\ell+1}}(\underline{\gamma}_{\ell+1},
\underline{\beta_{\ell+1}}|z)=\prod_{i=1}^{\ell+1}\prod_{j=1}^{\ell+1}
\Gamma(\imath \beta_{\ell+1,j}-\imath \gamma_{\ell+1,i}) e^{-\imath
z(\sum_{i=1}^{\ell+1}\gamma_{\ell+1,i}-\sum_{j=1}^{\ell+1}\beta_{\ell+1,j})},
\eqa acting on the space of functions of
$\underline{\gamma}=(\gamma_1,\ldots,\gamma_{\ell+1})$ as \bqa
\widehat{\CQ}^{\mathfrak{gl}_{\ell+1}}(z)\cdot
F(\underline{\gamma})= \int\limits_{{\cal
S}_{\ell+1}}\,\,\,\widehat{\CQ}^{\mathfrak{gl}_{\ell+1}}
(\underline{\gamma},\tilde{\underline{\gamma}}|z)
\,\,\,\,F(\tilde{\underline{\gamma}})\,\,\mu^{(\ell+1)}(\tilde{\underline{\gamma}})
\,\, \prod\limits_{ j=1}^{\ell+1}d\tilde{\gamma}_{j}\,. \eqa
\end{de}

\begin{prop} The $\mathfrak{gl}_{\ell+1}$-Whittaker
function satisfies  the following relation
\bqa\label{Baxterdual}
\widehat{\CQ}^{\mathfrak{gl}_{\ell+1}}(z)\cdot
\Psi^{\mathfrak{gl}_{\ell+1}}_{\underline{\gamma}_{\ell+1}}(\underline{x}_{\ell+1})
=e^{-e^{(x_{\ell+1,\ell+1}-z)}}
\Psi^{\mathfrak{gl}_{\ell+1}}_{\underline{\gamma}_{\ell+1}}
(\underline{x}_{\ell+1}). \eqa
\end{prop}
{\it Proof}.  We should prove that
\bqa\int\limits_{{\cal S}_{\ell+1}}\!\!\,\,e^{-\imath
z\sum_{i=1}^{\ell+1}(\lambda_{\ell+1,i}-\sum\gamma_{\ell+1,i})}
\,\,\cfrac
{\prod\limits_{i,j=1}^{\ell+1}\Gamma(\imath\gamma_{\ell+1,j}-\imath\lambda_{\ell+1,i})}
{\prod\limits_{i\neq j}
\Gamma(\imath\gamma_{\ell+1,j}-\imath\gamma_{\ell+1,i})}\,
\Psi^{\mathfrak{gl}_{\ell+1}}_{\underline{\gamma}_{\ell+1}}(\underline{x}_{\ell+1})\,
\,\,\prod_{j=1}^{\ell+1}d{\gamma}_{\ell+1,j}=\\
=\nonumber (2\pi)^{\ell+1}(\ell+1)!\,e^{-e^{(x_{\ell+1,\ell+1}-z)}}
\Psi^{\mathfrak{gl}_{\ell+1}}_{\underline{\lambda}_{\ell+1}}
(\underline{x}_{\ell+1}).\eqa
 Due to the orthogonality
condition (\ref{xorth}) this is equivalent to the following:
\bqa\int\limits_{\RR^{\ell+1}}\!\!\,e^{-e^{(x_{\ell+1,\ell+1}-z)}}
\overline{\Psi}^{\mathfrak{gl}_{\ell+1}}_
{\underline{\gamma}_{\ell+1}}(\underline{x}_{\ell+1})
\Psi^{\mathfrak{gl}_{\ell+1}}_{\underline{\lambda}_{\ell+1}}(\underline{x}_{\ell+1})
\,\prod_{j=1}^{\ell+1} d{x}_{\ell+1,j}\,=\\ \nonumber=e^{-\imath
z\sum_{i=1}^{\ell+1}(\lambda_{\ell+1,i}-\gamma_{\ell+1,i})}
\prod_{i=1}^{\ell+1}\prod_{j=1}^{\ell+1}
\Gamma(\imath\gamma_{\ell+1,j}-\imath\lambda_{\ell+1,i}).\eqa

Using the recursive relation (\ref{rec})  one can rewrite this as
\bqa\nonumber \int\limits_{\RR^{\ell+1}\times {\cal
S}_{\ell}\times{\cal
S}_{\ell}}\!\!\prod_{j=1}^{\ell+1}d{x}_{\ell+1,j}\,\prod_{j=1}^{\ell}d{\lambda}_{\ell,j}
\,\prod_{j=1}^{\ell}d{\gamma}_{\ell,j}\,\,\times\\ e^{-\imath
  x_{\ell+1,\ell+1}(\sum_{i=1}^{\ell+1}(\lambda_{\ell+1,i}-\gamma_{\ell+1,i})-
\sum_{k=1}^{\ell}(\lambda_{\ell,k}-\gamma_{\ell,k}))-e^{(x_{\ell+1,\ell+1}-z)}}\times
\\
\times  \cfrac{\prod\limits_{i=1}^{\ell+1}\prod\limits_{k=1}^{\ell}
\Gamma(\imath\la_{\ell,k}-\imath\lambda_{\ell+1,i})
\Gamma(\imath\gamma_{\ell+1,i}-\imath\gamma_{\ell,k})}
{(2\pi)^{2\ell}(\ell !)^2\prod\limits_{k\neq l}
\Gamma(\imath\lambda_{\ell,l}-\imath\lambda_{\ell,k})
\Gamma(\imath\gamma_{\ell,l}-\imath\gamma_{\ell,k})}\,
\overline{\Psi}^{\mathfrak{gl}_{\ell}}_{\underline{\gamma}_{\ell}}
(\underline{x}'_{\ell+1})
\Psi^{\mathfrak{gl}_{\ell}}_{\underline{\lambda}_{\ell}}(\underline{x}'_{\ell+1})
 \,=\\ \nonumber= e^{-\imath
z\sum_{i=1}^{\ell+1}(\lambda_{\ell+1,i}-\gamma_{\ell+1,i})}
\prod_{i=1}^{\ell+1}\prod_{j=1}^{\ell+1}
\Gamma(\imath\gamma_{\ell+1,j}-\imath\lambda_{\ell+1,i}), \eqa where
$\underline{x}'_{\ell+1}=(x_{\ell+1,1},\ldots,x_{\ell+1,\ell})$.
Using the orthogonality condition (\ref{xorth}) with respect to the
$\underline{x}'_{\ell+1}$ and integrating over
$\underline{\gamma}_{\ell}$ we see that (\ref{Baxterdual})
is equivalent to the following

\bqa\frac{1}{(2\pi)^\ell\,\ell\,!}
\int\limits_{-\infty}^{\infty}\!\!
dx_{\ell+1,\ell+1}\,\,e^{-\imath(x_{\ell+1,\ell+1}-z)\sum_{i=1}^{\ell+1}
(\lambda_{\ell+1,i}-\gamma_{\ell+1,i})-e^{x_{\ell+1,\ell+1}-z}}\cdot\\
\nonumber \int\limits_{{\cal
S}^{'}_{\ell}}\!\!\prod_{j=1}^{\ell}d{\lambda}_{\ell,j}\,\,
\cfrac{\prod\limits_{i=1}^{\ell+1}\prod\limits_{k=1}^{\ell}
\Gamma(\imath\la_{\ell,k}-\imath\lambda_{\ell+1,i})
\Gamma(\imath\gamma_{\ell+1,i}-\imath\la_{\ell,k})}
{\prod\limits_{k\neq l}
\Gamma(\imath\lambda_{\ell,l}-\imath\lambda_{\ell,k})}\,=\,
\prod_{i=1}^{\ell+1}\prod_{j=1}^{\ell+1}
\Gamma(\imath\gamma_{\ell+1,j}-\imath\lambda_{\ell+1,i}).\eqa
 Where the contour
of integration ${\cal S}^{'}_{\ell}$ in above formulas is deformed
so as to separate the sequences of poles going up $\{
\gamma_{\ell+1,j}+\imath k, j=1,\ldots,\ell+1, k=0,\ldots,\infty \}$
from the sequences of poles going down  $\{\lambda_{\ell+1,j}-\imath
k, j=1,\ldots,\ell+1, k=0,\ldots,\infty\}$ . We assume also that
$\gamma_{\ell+1,j}\neq\lambda_{\ell+1,k}$ for any  $j,k$. The last
identity is a simple consequence of  the following integral formula
due to Gustafson (see \cite{Gu}, Theorem 5.1, page 81):
$$
\frac{1}{(2\pi)^\ell}\int\limits_{{\cal
 S}^{'}_{\ell}}\!\!\prod_{j=1}^{\ell}d{\lambda_{\ell,j}}\,\,
\cfrac{\prod\limits_{i=1}^{\ell+1}\prod\limits_{k=1}^{\ell}
\Gamma(\imath\la_{\ell,k}-\imath\lambda_{\ell+1,i})
\Gamma(\imath\gamma_{\ell+1,i}-\imath \la_{\ell,k})}
{\prod\limits_{k\neq l}
\Gamma(\imath\lambda_{\ell,l}-\imath\lambda_{\ell,k})}\,=\,
\ell\,!\cfrac
{\prod\limits_{i=1}^{\ell+1}\prod\limits_{j=1}^{\ell+1}
\Gamma(\imath\gamma_{\ell+1,j}-\imath\lambda_{\ell+1,i})}
{\Gamma\Big(\sum\limits_{j=1}^{\ell+1}\imath
\gamma_{\ell+1,i}-\sum\limits_{i=1}^{\ell+1}\imath
\lambda_{\ell+1,i}\Big)}.
$$
 $\Box$

\begin{prop}\label{symrec}
The following symmetric recursive relation for
$\mathfrak{gl}_{\ell+1}$-Whittaker functions holds
\bqa\label{symrepOP}
\Psi^{\mathfrak{gl}_{\ell+1}}_{\underline{\gamma}_{\ell+1}}
(\underline{x}_{\ell+1})=e^{-\imath\gamma_{\ell+1,\ell+1}x_{\ell+1,\ell+1}}\,\,
\widehat{\CQ}^{\mathfrak{gl}_{\ell}}(x_{\ell+1,\ell+1})\cdot
\overline{\CQ^{\mathfrak{gl}_{\ell}}}(\gamma_{\ell+1,\ell+1})
\Psi^{\mathfrak{gl}_{\ell}}, \eqa where
$\underline{x}'_{\ell+1}=(x_{\ell+1,1},\ldots,x_{\ell+1,\ell})$,
$\underline{\gamma}'_{\ell+1}=(\gamma_{\ell+1,1},\ldots,\gamma_{\ell+1,\ell})$
and the action of the complex conjugated Baxter operator
$\overline{\CQ^{\mathfrak{gl}_{\ell}}}$ and its dual
$\widehat{\CQ}^{\mathfrak{gl}_{\ell}}$ is given by \bqa
\Big(\widehat{\CQ}^{\mathfrak{gl}_{\ell}}(x_{\ell+1,\ell+1})\cdot
\overline{\CQ^{\mathfrak{gl}_{\ell}}}(\gamma_{\ell+1,\ell+1})
\Psi^{\mathfrak{gl}_{\ell}}\Big)_{\underline{\gamma}_{\ell+1}}(\underline{x}_{\ell+1})=\\
=\int\prod_{j=1}^{\ell}
d{\gamma}_{\ell,j}\,\prod_{j=1}^{\ell}d{x}_{\ell,j}\,\,
\,\mu^{(\ell)}(\underline{\gamma}_{\ell})\,\widehat{\CQ}
^{\mathfrak{gl}_{\ell}}(\underline{\gamma}'_{\ell+1},
\underline{\gamma}_{\ell}|\,x_{\ell+1,\ell+1})\,
\overline{\CQ^{\mathfrak{gl}_{\ell}}}(\underline{x}'_{\ell+1},
\underline{x}_{\ell}|\,\gamma_{\ell+1,\ell+1}) \,
\Psi^{\mathfrak{gl}_{\ell}}_{\underline{\gamma}_{\ell}}(\underline{x}_{\ell}).
\nonumber\eqa
\end{prop}
{\it Proof}. Let us start with the Mellin-Barnes recursive relation
\be\nonumber
\Psi^{\mathfrak{gl}_{\ell+1}}_{\underline{\gamma}_{\ell+1}}
(\underline{x}_{\ell+1})=\\\int
\,\prod_{j=1}^{\ell}d{\gamma}_{\ell,j}\,\,\mu^{(\ell)}
(\underline{\gamma}_{\ell})\,\,
\prod_{i=1}^{\ell+1}\prod_{j=1}^{\ell}\,\Gamma(\imath
\gamma_{\ell,j} -\imath \gamma_{\ell+1,i}) e^{-\imath
x_{\ell+1,\ell+1}
(\sum_{i=1}^{\ell+1}\gamma_{\ell+1,i}-\sum_{i=1}^{\ell}\gamma_{\ell,i})}
\Psi^{\mathfrak{gl}_{\ell}}_{\underline{\gamma}_{\ell}}(\underline{x}'_{\ell+1}).
\ee Using the properties of the Baxter operator we have
$$
\Psi^{\mathfrak{gl}_{\ell+1}}_{\underline{\gamma}_{\ell+1}}(\underline{x}_{\ell+1})
=e^{-\imath \gamma_{\ell+1,\ell+1}x_{\ell+1,\ell+1}}
 \int \prod_{j=1}^{\ell}d{x}_{\ell,j}\,\,
 \overline{\CQ^{\mathfrak{gl}_{\ell}}}
(\underline{x}'_{\ell+1},\underline{x}_{\ell}|\gamma_{\ell+1,\ell+1})\times
$$
$$
\times \Big(\int\prod_{j=1}^{\ell} d{\gamma}_{\ell,j}\,\,
\mu^{(\ell)}(\underline{\gamma}_{\ell})\,\,
\prod_{i=1}^{\ell}\prod_{j=1}^{\ell}\,\Gamma(\imath \gamma_{\ell,j}-
\imath \gamma_{\ell+1,i}) e^{-\imath
x_{\ell+1,\ell+1}(\sum_{i=1}^{\ell}\gamma_{\ell+1,i}-\sum_{i=1}^{\ell}\gamma_{\ell,i})}\,
\Psi^{\mathfrak{gl}_{\ell}}_{\underline{\gamma}_{\ell}}(\underline{x}_{\ell})\Big)=$$
$$
=e^{-\imath\gamma_{\ell+1,\ell+1}x_{\ell+1,\ell+1}}\,\,\widehat{\CQ}
^{\mathfrak{gl}_{\ell}}(x_{\ell+1,\ell+1})\cdot
\overline{\CQ^{\mathfrak{gl}_{\ell}}}(\gamma_{\ell+1,\ell+1})
\Psi^{\mathfrak{gl}_{\ell}}.
$$
Note that one can equally start with a Givental recursive relation
and use the eigenvalue property (\ref{Baxterdual})
of the dual Baxter operator. $\Box$

The Givental and Mellin-Barnes recursions are easily obtained from  the
symmetric recursion  (\ref{symrepOP}).  This
provides a direct and inverse transformation of the  Givental
representation into the Mellin-Barnes one. Moreover, this leads to
a family  of the intermediate Givental-Mellin-Barnes representations.
Indeed, to obtain
 $\Psi^{\mathfrak{gl}_{\ell+1}}_{\underline{\gamma}_{\ell+1}}(\underline{x}_{\ell+1})$
from $\Psi^{\mathfrak{gl}_{\ell}}_{\underline{\gamma}_{\ell}}(\underline{x}_{\ell})$
one can either use the integral operator $Q^{\mathfrak{gl}_{\ell+1}}_{\mathfrak{gl}_{\ell}}
(\underline{x}_{\ell+1},\underline{x}_{\ell}|\gamma_{\ell+1,\ell+1})$
or the  integral operator $
\widehat{Q}^{\mathfrak{gl}_{\ell+1}}_{\mathfrak{gl}_{\ell}}
(\underline{\gamma}_{\ell+1},\underline{\gamma}_{\ell}|x_{\ell+1,\ell+1})$.
 This leads to the following family
of mixed Mellin-Barnes-Givental integral representations  of
$\mathfrak{gl}_{\ell+1}$-Whittaker function \bqa\label{interpol}
\Psi^{\mathfrak{gl}_{\ell+1}}=Q^{(\epsilon_1)}\cdot
Q^{(\epsilon_2)}\cdots
  Q^{(\epsilon_{\ell})}\,\, \Psi^{\mathfrak{gl}_{1}},\qquad \epsilon=L,R\,\,,
\eqa
where $Q^{(L)}$ is the integral operator with the integral kernel
$Q^{\mathfrak{gl}_{k+1}}_{\mathfrak{gl}_{k}}$,
$Q^{(R)}$ is the integral operator with the integral kernel
$\widehat{Q}^{\mathfrak{gl}_{k+1}}_{\mathfrak{gl}_{k}}$
 and the integral operators act on $\underline{\gamma}$- or
$\underline{x}$-variables depending on $\epsilon_i$. Various
 choices of $\{\epsilon_i\}$ in (\ref{interpol}) provide
various  integral representations of
$\mathfrak{gl}_{\ell+1}$-Whittaker  function.

\section{ Archimedean factors in  Rankin-Selberg method}

In this section we apply  the  dual recursion operator and Baxter
operators discussed in the previous section to simplify
calculations  of the correction factors arising in the
Rankin-Selberg method applied to $GL(\ell+1)\times GL(\ell+1)$ and
$GL(\ell+1)\times GL(\ell)$.
 Note that these calculations  are  an important step in the proof
of the functional equations for the corresponding
 automorphic $L$-functions using the Rankin-Selberg approach.
Explicit expressions for these correction factors in terms of
Gamma-functions were conjectured by Friedberg-Bump and Bump and
proved later by Stade \cite{St1},\, \cite{St2}. The proofs in
\cite{St1},\, \cite{St2} are  based on a recursive generalization of
the  integral representation of $\mathfrak{gl}_{\ell+1}$-Whittaker
functions, $\ell=2$ first derived by Vinogradov and Takhtadzhyn
\cite{VT}. The recursion in \cite{St1},\,\cite{St2} changes the rank
by two $\ell-1\to \ell+1$. It was noted in \cite{GKLO} that this
recursion is basically the Givental recursion applied twice.

In this section we will demonstrate that
using the recursive properties of the Mellin-Barnes representation and
the dual Baxter operator  one can give a one-line proof of  Friedberg-Bump and
Bump conjectures. We start with a brief description of the relevant
facts about automorphic $L$-functions, Rankin-Selberg method
and Bump-Freidberg and Bump conjectures. For more details see e.g. \cite{Bu}, \cite{Go}.

Let $\mathbb{A}$ be the adele ring of  $\mathbb{Q}$ and $G$ be a
reductive Lie group. An automorphic representation $\pi$ of
$G(\mathbb{A})$ can be characterized by an automorphic form
$\phi_{\pi}$ such that it is an eigenfunction of  any element of the
global Hecke algebra $\CH(G(\mathbb{A}))$. The global Hecke algebra can
be represented as a  product $\CH(G(\mathbb{A}))=(\otimes_p
\CH_p)\otimes \CH_{\infty}$ of  the local non-Archimedean Hecke
algebras $\CH_p=\CH(G(\mathbb{Q}_p), G(\mathbb{Z}_p))$ for each
prime $p$ and an Archimedean Hecke algebra
$\CH_{\infty}=\CH(G(\mathbb{R}),K)$ where $K$ is a maximal compact
subgroup in $G(\RR)$. Local Hecke algebra $\CH_p$
 is isomorphic to a representation ring of  a simply connected
 complex Lie group ${}^LG_0$, Langlands dual to $G$
(e.g.  $A_{\ell}$, $B_{\ell}$,  $C_{\ell}$,  $D_{\ell}$ are dual to
 $A_{\ell}$, $C_{\ell}$,  $B_{\ell}$,  $D_{\ell}$ respectively).
For each unramified representation  of $G(\mathbb{Q}_p)$ one can
define an action of $\CH_p$ such that an automorphic form
$\phi_{\pi}$  is a common eigenfunction of all elements of  $\CH_p$
for all  primes $p$ and thus defines  a set of homomorphisms
$\CH_p\to \mathbb{C}$. Identifying local Hecke algebras with the
representation ring of ${}^LG_0$ one can describe this set of
homomorphisms as a  set  of conjugacy classes $g_p$  in ${}^LG_0$.

Given a finite-dimensional representation
$\rho_V: {}^LG_0\to GL(V,\CC)$ one can
construct an $L$-function corresponding to an automorphic form
$\phi$ in the form of the Euler product as follows
\bqa\label{autoL}
L(s,\phi,\rho_V)=\mathop{{\prod}'}_{p}
L_p(s,\phi,\rho_V)=\mathop{{\prod}'}_{p} \,
\det_{V} (1-\rho_V(g_p)\, p^{-s})^{-1},
\eqa
where $\prod_p'$ is a product over primes $p$ such that
the corresponding representation of $G(\mathbb{Q}_p)$ is not ramified.
It is natural to complete the product by including  local
$L$-factors corresponding to Archimedean and ramified places.
$L$-factors for  ramified representations can be taken trivial.
For the Archimedean place the
Hecke eigenfunction property is  usually replaced  by the eigenfunction property
with respect to  the ring of invariant
differential operators on $G(\mathbb{R})$. The corresponding
eigenvalues are described by a conjugacy class  $t_{\infty}$ in the Lie
 algebra ${}^L\mathfrak{g}_0={\rm Lie}({}^LG_0)$. The Archimedean
 $L$-factor  is given by \cite{Se}
\bqa\label{archL}
L_{\infty}(s,\phi,\rho_V)=\prod_{j=1}^{\ell+1}
\Big(\pi^{-\frac{s-\alpha_j}{2}}\,
\Gamma\Big(\frac{s-\alpha_j}{2}\Big)\Big)
=\det_V\Big(\pi^{-\frac{s-\rho_V(t_{\infty})}{2}}\,
\Gamma\Big(\frac{s-\rho_V(t_{\infty})}{2}\Big)\Big),
\eqa
where $\rho_V(t_{\infty})={\rm diag}(\alpha_1,\ldots \alpha_{\ell+1})$.
The complete $L$-function
\bqa
\Lambda(s,\phi,\rho)=L(s,\phi,\rho)L_{\infty}(s,\phi,\rho),
\eqa
should satisfy the  functional equation of the form
$$
\Lambda(1-s,\phi,\rho)=\epsilon(s,\phi,\rho)
\Lambda(s,\phi_{\pi^{\vee}},\rho^{\vee}),
$$
where $\epsilon$-factor is of the exponential form
$\epsilon(s,\phi,\rho)=A\,B^s$ and $\pi^{\vee}$, $\rho^{\vee}$
are dual to $\pi$, $\rho$.

In Rankin-Selberg method one considers  automorphic $L$-functions
associated with automorphic representations of the
products $G\times \tilde{G}$ of reductive groups.
Let  $\rho_V: {}^LG_0\to
{\rm End}(V)$,  $\tilde{\rho}_{\tilde{V}}: {}^L\tilde{G}_0\to {\rm
  End}(\tilde{V})$ be finite-dimensional representations of  dual groups
and let $g_p\in {}^LG_0$, $\tilde{g}_p\in {}^L\tilde{G}_0$ be representatives of the
conjugacy classes corresponding to automorphic forms $\phi$ and
$\tilde{\phi}$. One defines  $L$-function
 $L(s,\pi\times \tilde{\pi},\rho\times \tilde{\rho})$ as follows
\bqa\label{RSLfunc}
L(s,\phi\times \tilde{\phi}, \rho\times \tilde{\rho})=\mathop{{\prod}'}_{p}
\,\det_{V\otimes \tilde{V}} (1-\rho_{V}(g_p)\otimes \tilde{\rho}_V(\tilde{g}_p)
\, p^{-s})^{-1}.
\eqa
$L$-function (\ref{RSLfunc}) up to a  correction factor
can be naturally written as  an integral of the product of  automorphic forms
$\phi$  and $\tilde{\phi}$  with a simple kernel function.
Given an explicit expression for the correction factor,
this integral representation can be an important tool for  studying
analytic properties of $L(s,\phi\times \tilde{\phi})$ as a function of
$s$.

In the following we consider Rankin-Selberg method in case of
$G\times\widetilde{G}$ being either $\qquad$ $GL(\ell+1)\times GL(\ell+1)$
or  $GL(\ell+1)\times GL(\ell)$ with $\rho$ and
$\tilde{\rho}$  being standard representations. We
start with the case of $GL(\ell+1)\times GL(\ell+1)$.
Consider the following zeta-integral
\bqa\label{onetwo}
Z(s,\phi \times \tilde{\phi})=\int_{GL(\ell+1,\mathbb{Q})Z^{(\ell+1)}
_{\mathbb{A}}\backslash
 GL(\ell+1,\mathbb{A})}
\phi(g)\tilde{\phi}(g) \mathcal{E}(g,s)\,dg,
\eqa
where the Eisenstein series is
\bqa
\mathcal{E}(g,s)=\zeta((\ell+1)s)\,
\sum_{\gamma \in P(\ell+1,\ell,
\mathbb{Z})\backslash GL(\ell+1,\mathbb{Z})}\,f_s(\gamma g).
\eqa
Here  $Z^{(\ell+1)}_{\mathbb{A}}$ is the  center of $GL(\ell+1,\mathbb{A})$,
$\zeta(s)=\sum_{n=1}^{\infty}n^{-s}$ is the Riemann zeta-function
 and
$$f_s\in {\rm Ind}^{GL(\ell+1,\mathbb{A})}_{P(\ell+1,\ell,\mathbb{A})}\,\,\delta_P^s,$$
where  $\delta_P$ denotes the modular function of the parabolic subgroup
$P(\ell+1,\ell,\mathbb{A})$ of $GL(\ell+1,\mathbb{A})$
with the Levi factor $GL(\ell,\mathbb{A})\times GL(1,\mathbb{A})$.

Using the Rankin-Selberg unfolding technique
 (\ref{onetwo})  can be represented in the form
$$
Z(s,\phi\times \tilde{\phi})=L(s,\phi \times \tilde{\phi})
\Psi(s,\phi\times \tilde{\phi}),
$$
where the correction factor $\Psi(s,\phi\times \tilde{\phi})$ is  a
convolution  of two $\mathfrak{gl}_{\ell+1}$-Whittaker functions.
The Bump-Freidberg  conjecture proved in  \cite{St1} claims that
$\Psi(s,\phi\times \tilde{\phi})$ is equal
to the  Archimedean local $L$-factor.

\begin{te}[Bump-Freidberg-Stade]
\label{BFS}
\bqa\label{BumpStone}
\Psi(s, \phi\times \tilde{\phi})=L_{\infty}(s, \phi\times
\tilde{\phi})=\prod_{j=1}^{\ell+1}\prod_{k=1}^{\ell+1}\,\pi^{-\frac{s-\alpha_j-
\tilde{\alpha}_k}{2}}
\Gamma\Big(\frac{s-\alpha_j-\tilde{\alpha}_k}{2}\Big),
\eqa
where $\rho_V(t_{\infty})={\rm diag}(\alpha_1,\ldots \alpha_{\ell+1})$ and
$\tilde{\rho}_V(\tilde{t}_{\infty})={\rm diag}(\tilde{\alpha}_1,\ldots
\tilde{\alpha}_{\ell+1})$ correspond to the automorphic
representations $\phi$ and $\tilde{\phi}$ as in  (\ref{archL}).
\end{te}

The proof of Theorem can be reduced to the following identity
proved by Stade (we rewrite Theorem 1.1, \cite{St2}
in our notations).

\begin{lem} The following integral relation holds
\bqa
\int\limits_{\RR^{\ell+1}}\!\prod_{j=1}^{\ell+1}d{x}_{\ell+1,j}\,\,
e^{-e^{x_{\ell+1,\ell+1}}}
\overline{\Psi^{\mathfrak{gl}_{\ell+1}}}_{\underline{\gamma}_{\ell+1}}
(\underline{x}_{\ell+1})
\Psi^{\mathfrak{gl}_{\ell+1}}_{\underline{\la}_{\ell+1}+\underline{t}}
(\underline{x}_{\ell+1})
=\prod_{k=1}^{\ell+1}\prod_{j=1}^{\ell+1}\Gamma(\imath t+\imath
\lambda_{\ell+1,k}-\imath \gamma_{\ell+1,j}), \eqa where
$\underline{t}=(t,\ldots ,t)\in \mathbb{R}^{\ell+1}$.
\end{lem}

{\it Proof}. The proof readily follows from the proof of Proposition 3.2.
 $\Box$.

Next we consider the  Rankin-Selberg method for
$GL(\ell+1)\times GL(\ell)$, $\rho$ and
$\tilde{\rho}$  being standard representations of
$GL(\ell+1)$ and $GL(\ell)$. In this case one has to study
the following integral
\bqa\label{intrepRS}
Z(s,\phi \times \tilde{\phi})=\int_{GL(\ell,\ZZ)Z^{(\ell)}_{\mathbb{A}}\backslash
  GL(\ell,\mathbb{A})}\phi(\begin{pmatrix}
  g& \\
& 1\end{pmatrix})\,\tilde{\phi}(g)\, |{\rm det}(g)|^{s-1/2} dg,
\eqa
where $Z^{(\ell)}$ is the  center of $GL(\ell,\mathbb{A})$.
Using the Rankin-Selberg unfolding technique, the integral
 (\ref{intrepRS})  can be represented in the form
$$
Z(s,\phi\times \tilde{\phi})=L(s,\phi \times \tilde{\phi})
\Psi(s,\phi\times \tilde{\phi}),
$$
where the correction factor $\Psi(s,\phi\times \tilde{\phi})$ is a
convolution  of $\mathfrak{gl}_{\ell+1}$-  and
$\mathfrak{gl}_{\ell}$-Whittaker functions.
The Bump conjecture proved in  \cite{St1} claims that
$\Psi(s,\phi\times \tilde{\phi})$ is equal
to the  Archimedean local $L$-factor.
\begin{te}[Bump-Stade]
\label{BS}
\bqa\label{BumpStoneone}
\Psi(s, \phi\times \tilde{\phi})=L_{\infty}(s, \phi\times
\tilde{\phi})=\prod_{j=1}^{\ell+1}\prod_{k=1}^{\ell}\,\pi^{-\frac{s-\alpha_j-
\tilde{\alpha}_k}{2}}
\Gamma\Big(\frac{s-\alpha_j-\tilde{\alpha}_k}{2}\Big),
\eqa
where $\rho_V(t_{\infty})={\rm diag}(\alpha_1,\ldots \alpha_{\ell+1})$ and
$\tilde{\rho}_V(\tilde{t}_{\infty})={\rm diag}(\tilde{\alpha}_1,\ldots
\tilde{\alpha}_{\ell})$ correspond to the automorphic
representations $\phi$ and $\tilde{\phi}$ as in  (\ref{archL}).
\end{te}

The proof of the theorem is equivalent to the proof of the following
integral identity ( we rewrite the Theorem 3.4,
\cite{St2} using our notations).

\begin{lem}
\bqa\label{orthGL}\int\limits_{\RR^{\ell+1}}
\prod_{j=1}^{\ell+1}d{x}_{\ell+1,j}\,\,
\overline{\Psi}^{\mathfrak{gl}_{\ell}}_{\underline{\gamma}_{\ell}}
(\underline{x}'_{\ell+1})
\Psi^{\mathfrak{gl}_{\ell+1}}_{\underline{\la}_{\ell+1}+\underline{t}}
(\underline{x}_{\ell+1})
=\\
=\nonumber
\delta\Big(\,\i(\ell+1)t+\i\sum_{i=1}^{\ell+1}\la_{\ell+1,i}\,-\,
\i\sum_{k=1}^{\ell}\gamma_{\ell,k}\Big)\,
\prod_{i=1}^{\ell+1}\prod_{k=1}^{\ell} \Gamma(\i
t+\i\la_{\ell+1,i}-\i\gamma_{\ell,k})\eqa where $\underline{t}={\rm
diag}(t,\ldots,t)\in\mathbb{R}^{\ell+1}$,
$\underline{x}'_{\ell+1}=(x_{\ell+1,1},\ldots,x_{\ell+1,\ell})$ and
$\delta(x)$ is the Dirac $\delta$-function.
\end{lem}
{\it Proof}.  To verify this statement we substitute into the l.h.s.
of (\ref{orthGL}) the following recursive relation
$$
\Psi^{\mathfrak{gl}_{\ell+1}}_{\underline{\la}_{\ell+1}+\underline{t}}
(\underline{x}_{\ell+1})=\widehat{Q}(x_{\ell+1,\,\ell+1})\cdot
\Psi^{\mathfrak{gl}_{\ell}}_{\underline{\la}_{\ell}}
(\underline{x}'_{\ell})
$$
Then applying orthogonality relation from Theorem 2.1 and
integrating over $x_{\ell+1,\,\ell+1}$ we obtain the r.h.s. (\ref{orthGL}) $\Box$

Let us stress that one should not expect to have expressions
for  $\Psi(s,\phi\times \tilde{\phi})$ as products of Gamma-functions
for more general cases $GL(\ell+n)\times GL(\ell)$, $n>1$.
 From the point of view of Mellin-Barnes recursive construction,
$\Psi(s,\phi\times \tilde{\phi})$  are the kernels of recursive operators
corresponding to the change of rank $\ell\to \ell+n$ and
thus are given by  compositions of elementary recursive operators.
This leads to  general  expressions for $\Psi(s,\phi\times \tilde{\phi})$
in terms of the  integrals of the products of Gamma-functions.
Let us remark that in this paper we consider Rankin-Selberg method
as a method for studying  properties of  matrix elements
of the natural (recursive) operators acting in the space of automorphic forms.
One can expect that this point of view might be useful in the
investigation  of  other properties of automorphic $L$-functions.

Let us comment on Stade  proof of
Theorems \ref{BFS}, \ref{BS}. The proof in \cite{St1}, \cite{St2}
is based on the recursive relation connecting $\mathfrak{gl}_{\ell+1}$- and
$\mathfrak{gl}_{\ell-1}$-Whittaker functions. Below we derive this recursion
from  the following form of the  Givental recursion.

\begin{prop} The following recursive relations for
$\mathfrak{gl}_{\ell+1}$-Whittaker functions holds:
 \bqa\label{Gtwice}
\Psi^{\mathfrak{gl}_{\ell+1}}_{\lambda_1,\ldots,\lambda_{\ell+1}}
(\underline{x}_{\ell+1})= \int_{\RR^{\ell-1}}
\prod_{i=1}^{\ell-1}dx_{\ell-1,i}\,\,
Q^{\mathfrak{gl}_{\ell+1}}_{\mathfrak{gl}_{\ell-1}}
(\underline{x}_{\ell+1},\underline{x}_{\ell-1}|\lambda_{\ell+1},\lambda_{\ell})
\Psi_{\lambda_1,\ldots,\lambda_{\ell-1}}^{\mathfrak{gl}_{\ell-1}}
(\underline{x}_{\ell-1}),\eqa \bqa\nonumber
Q^{\mathfrak{gl}_{\ell+1}}_{\mathfrak{gl}_{\ell-1}}
(\underline{x}_{\ell+1},\underline{x}_{\ell-1}|\lambda_{\ell+1},\lambda_{\ell})=\eqa
\bqa \label{GtwiceQ}
=\int_{\RR^{\ell}}\,\prod_{j=1}^{\ell}d{x}_{\ell,j}\,
\exp\Big\{\i\la_{\ell+1}\Big(\,\sum_{i=1}^{\ell+1}x_{\ell+1,i}-
\sum_{k=1}^{\ell}x_{\ell,k}\Big)-
\sum_{k=1}^{\ell}\bigl(e^{x_{\ell+1,k}-x_{\ell,k}}+
e^{x_{\ell,k}-x_{\ell+1,k+1}}\bigr)\,+\\ \nonumber
+\i\la_\ell\Big(\,\sum_{k=1}^{\ell}x_{\ell,k}-
\sum_{j=1}^{\ell-1}x_{\ell-1,j}\Big)-
\sum_{k=1}^{\ell-1}\bigl(e^{x_{\ell,k}-x_{\ell-1,k}}+
e^{x_{\ell-1,k}-x_{\ell,k+1}}\bigr)\,\Big\}.\eqa
\end{prop}
{\it Proof}. The recursive relation (\ref{Gtwice}) is the Givental
recursive relation (\ref{recursgl})  applied twice $\Box$

\begin{te}[Stade]
The following recursion relation for
$\mathfrak{gl}_{\ell+1}$-Whittaker functions holds
\bqa\label{recstad}
\Psi^{\mathfrak{gl}_{\ell+1}}_{\lambda_1,\ldots,\lambda_{\ell+1}}
(\underline{x}_{\ell+1})=\int_{\mathbb{R}^{\ell-1}}
\prod_{j=1}^{\ell-1}d{x}_{\ell-1,j}\,\, K_{\ell+1,\,\ell-1}
(\underline{x}_{\ell+1},\,\underline{x}_{\ell-1}|\,\lambda_{\ell+1},\lambda_{\ell})
\,\Psi^{\mathfrak{gl}_{\ell-1}}_{\lambda_1,\ldots\lambda_{\ell-2}}
(\underline{x}_{\ell-1}),
 \eqa where
 $K_{\ell+1,\,\ell-1}
(\underline{x}_{\ell+1},\,\underline{x}_{\ell-1}|\,\underline{\lambda})$
 is given by the following explicit formula.
\bqa\label{kerstade} K_{\ell+1,\,\ell-1}
(\underline{x_{\ell+1}},\,\underline{x_{\ell-1}}|\,\lambda)=
2^{1-\ell}\exp\Big\{\,\frac{\i(\la_{\ell}+\la_{\ell+1})}{2}\Big(
\sum_{i=1}^{\ell+1}x_{\ell+1,i}-
\sum_{j=1}^{\ell-1}x_{\ell-1,j}\Big)\,\Big\}\times \\
\nonumber\times  \prod_{i=1}^\ell\,K_{\i(\la_\ell-\la_{\ell+1})}\,
\Big(\,2\sqrt{\bigl(e^{x_{\ell+1,i}}+e^{x_{\ell-1,i-1}}\bigr)
\bigl(e^{-x_{\ell+1,i+1}}+e^{-x_{\ell-1,i}}\bigr)}\,\Big)\eqa Here
we use the following  integral representation for the Macdonald
function
$$
K_\nu(y)=\int_0^\infty\frac{dt}{t}\,\,t^\nu e^{-y(t+t^{-1})/2}
$$
\end{te}

{\it Proof}: At first we substitute into the expression for
$K_{\ell+1,\,\ell-1}$ the  integral representation with integration
variables $t_{i}$ for Macdonald functions
$K_{\i(\la_\ell-\la_{\ell+1})}$. Then we make the following change
of variables $t_i$. \bqa
t_1=e^{x_{\ell,1}}\sqrt{\frac{e^{-x_{\ell+1,2}}+e^{-x_{\ell-1,1}}}{e^{x_{\ell+1,1}}}}
\eqa \bqa\nonumber
t_k=e^{x_{\ell,k}}\sqrt{\frac{e^{-x_{\ell+1,k+1}}+e^{-x_{\ell-1,k}}}
{e^{x_{\ell+1,k}}\,+\,e^{x_{\ell-1,k-1}}}}\hspace{1cm}
t_{\ell}=e^{x_{\ell,\ell}}\sqrt{\frac{e^{-x_{\ell+1,\ell+1}}}
{e^{x_{\ell+1,\ell}}+e^{x_{\ell-1,\ell-1}}}},\eqa for
$k=1,\ldots,\ell$ and $j=1,\ldots,\ell-1$. Thus we obtain the
following identity between the kernels
 \bqa
K_{\ell+1,\,\ell-1}
(\underline{x}_{\ell+1},\,\underline{x}_{\ell-1}|\,\lambda)\,=\,
Q^{\mathfrak{gl}_{\ell+1}}_{\mathfrak{gl}_{\ell-1}}
(\underline{x}_{\ell+1},\,\underline{x}_{\ell-1}|\,\lambda)\eqa This
reduces Stade recursion to the Givental recursive procedure  $\Box$

The appearance of  the Gamma-functions both in
the Mellin-Barnes integral representation of the
$\mathfrak{gl}_{\ell+1}$-Whittaker functions and
in the expressions for the Archimedean $L$-factors is not accidental.
In the next section we explain this connection by relating
the constructed Baxter operator with a universal Baxter operator
considered as an element of the  Archimedean Hecke algebras $\CH(G(\mathbb{R}),K)$
where $K$ is a maximal compact subgroup of $G(\mathbb{R})$.

\section{ Universal Baxter operator}

\subsection{Universal Baxter operator in $\CH(G(\mathbb{R}),K)$}

In this section we will argue  the Baxter $\mathcal{Q}$-operator
for the  $\mathfrak{gl}_{\ell+1}$-Toda chain can  (and should) be
considered as a realization of
the universal Baxter operator considered as elements of the
spherical  Hecke algebra $\CH(GL(\ell+1,\mathbb{R}),K)$, $K$ being a
maximal compact subgroup of $GL(\ell+1,\mathbb{R})$.
We  also consider   non-Archimedean analogs
of the universal Baxter operator as an element of
a local Hecke algebra $\mathcal{H}(GL(\ell+1,\mathbb{Q}_p),
GL(\ell+1,\mathbb{Z}_p))$.
Both in   Archimedean and  non-Archimedean cases the eigenvalues of
the Baxter $\mathcal{Q}$-operators
acting on $\mathfrak{gl}_{\ell+1}$-Whittaker functions
are  given by the corresponding local $L$-factors.

Let us start with the definition  of the
spherical  Hecke algebra
$\mathcal{H}_{\infty}=\mathcal{H}(G(\mathbb{R}),K)$, where $K$ is a
maximal compact subgroup of $G(\mathbb{R})$.
 Algebra $\mathcal{H}_{\infty}$ is defined
as an algebra of $K$-biinvariant functions
 on $G$,  $\phi(g)=\phi(k_1gk_2)$, $k_1,k_2\in K$  acting by a convolution
\bqa \phi*f(g)=\int_G \phi(g\tilde{g}^{-1})\,f(\tilde{g})
d\tilde{g}. \eqa To ensure the convergence of the integrals one
usually imposes the condition of compact support on $K$-biinvariant
functions. We will consider slightly more general class of
exponentially decaying functions.\footnote{ This should  be compared
with  the use of  exponentially decreasing functions instead of
functions with compact support in the Mathai-Quillen construction of
the representative of the Thom class.}.

By the multiplicity one theorem \cite{Sha}, there is a unique
smooth  spherical vector $\<k|$ in a principal series irreducible
representation $\mathcal{V}_{\underline{\gamma}}={\rm
Ind}_{B_-}^G\,\chi_{\underline{\gamma}}$. The action of a
$K$-biinvariant function $\phi$ on the spherical vector $\<k|$ in
$\mathcal{V}_{\underline{\gamma}}$ is given by  the multiplication
by  a character $\Lambda_{\phi}$ of the Hecke algebra: \be
\phi*\<k|\equiv\int_G dg
\phi(g^{-1})\,\<k|\pi_{\gamma}(g)=\Lambda_{\phi}(\gamma)\<k|. \ee In
particular, the  elements $\phi$ of the Hecke algebra should act by
convolution on the Whittaker function  as follows
\bqa\label{phiEigen}
\phi*\Phi^{\mathfrak{gl}_{\ell+1}}_{\underline{\gamma}}
(g)=\Lambda_{\phi}(\underline{\gamma})
 \Phi^{\mathfrak{gl}_{\ell+1}}_{\underline{\gamma}}
(g),\qquad \phi\in \CH_{\infty}. \eqa Here the Whittaker function
$\Phi_{\underline{\gamma}}^{\mathfrak{gl}_{\ell+1}}$ is considered
as a function on $G$ such that \be\label{equivar}
\Phi^{\mathfrak{gl}_{\ell+1}}_{\underline{\gamma}}(kan)=\chi_{N_{-}}(n)\,
 \Phi^{\mathfrak{gl}_{\ell+1}}_{\underline{\gamma}}(a),
\ee where $kan\in KAN_-\to G$ is the Iwasawa decomposition.

In the previous section we construct the Baxter integral operator
acting on the $\mathfrak{gl}_{\ell+1}$- Whittaker function
(considered as a function on the subspace $A$ of the diagonal
matrices)  as \be\label{BAxterred}
\CQ^{\mathfrak{gl}_{\ell+1}}(\lambda)\,
\cdot\Psi^{\mathfrak{gl}_{\ell+1}}_{\underline{\gamma}}
(\underline{x})=\prod_{j=1}^{\ell+1}\,
\pi^{-\frac{\imath\lambda-\imath\gamma_j}{2}}\,
\Gamma\Big(\frac{\imath \lambda-\imath \gamma_{j}}{2}\Big)
\,\,\Psi^{\mathfrak{gl}_{\ell+1}}_{\underline{\gamma}}(\underline{x}),
\ee where the kernel of  the operator
$\CQ^{\mathfrak{gl}_{\ell+1}}(\lambda)$ is given by
$$\CQ^{\mathfrak{gl}_{\ell+1}}
(\underline{x},\underline{y}|\l)=2^{\ell+1}\exp\Big\{\imath\l\sum_{i=1}^{\ell+1}
(x_i-y_i)-\pi
\sum_{k=1}^{\ell}\Big(e^{2(x_k-y_k)}+e^{2(y_k-x_{k+1})}\Big)- \pi
e^{2(x_{\ell+1}-y_{\ell+1})}\Big\}$$ Note that here we use a
parametrization of  Baxter operator naturally arising in the
description of Whittaker functions in terms of Iwasawa
decomposition.  In this section we will use only this type of the
parametrization and drop the tildes in the corresponding notations
(see \ref{SpherW} and Lemma \ref{Iwarep}). We also take coupling
constants in Toda chain $g_i=\pi$ to agree  with the standard
normalizations in  Representation theory.

Let us recall that we introduce $\mathfrak{gl}_{\ell+1}$-Whittaker
function
$\Psi^{\mathfrak{gl}_{\ell+1}}_{\underline{\gamma}}(\underline{x})$
as a matrix element multiplied by the factor $\exp(-\<\rho,x\>)$
(see (\ref{pairgl}), (\ref{SpherW})). In the construction of the
universal Baxter operator it is more natural to consider a modified
Whittaker functions $\Phi^{\mathfrak{gl}_{\ell+1}}$ equal to  the
matrix elements itself \be
\Phi^{\mathfrak{gl}_{\ell+1}}_{\underline{\gamma}}(\underline{x})=
e^{\<\rho,x\>}\,
\Psi^{\mathfrak{gl}_{\ell+1}}_{\underline{\gamma}}(\underline{x}).
\ee Define a modified Baxter $\CQ$-operator:
$$\CQ_0^{\mathfrak{gl}_{\ell+1}}(\lambda)
=e^{\<\rho,{x}\>}\CQ^{\mathfrak{gl}_{\ell+1}}(\lambda)e^{-\<\rho,{x}\>}.$$
It has the kernel \be \CQ_0^{\mathfrak{gl}_{\ell+1}}
(\underline{x},\underline{y}|\l)=2^{\ell+1}\exp\Big\{\sum_{j=1}^{\ell+1}(\i
\l +\rho_j) (x_j-y_j)-\\ \nonumber -\pi
\sum_{k=1}^{\ell}\Big(e^{2(x_k-y_k)}+e^{2(y_k-x_{k+1})}\Big)- \pi
e^{2(x_{\ell+1}-y_{\ell+1})}\Big\}, \ee
 where $\rho\in \RR^{\ell+1}$ , with $\rho_j=\frac{\ell}{2}+1-j,\,\,\,j=1,\ldots,\ell+1$,
 and it acts on the modified
Whittaker functions as follows \be\label{BAxterredone}
\CQ_0^{\mathfrak{gl}_{\ell+1}}(\lambda)\,
\cdot\Phi^{\mathfrak{gl}_{\ell+1}}_{\underline{\gamma}}
(\underline{x})=\prod_{j=1}^{\ell+1}\,
\pi^{-\frac{\imath\lambda-\imath\gamma_j+\rho_j}{2}}\,
\Gamma\Big(\frac{\imath \lambda-\imath \gamma_{j}+\rho_j}{2}\Big)
\,\,\Phi^{\mathfrak{gl}_{\ell+1}}_{\underline{\gamma}}(\underline{x}).
\ee We would like to find an element $\phi_{\CQ_0(\l)}$ in
$\mathcal{H}_{\infty}$  such that the following relation holds \be
\phi_{\CQ_0(\lambda)}\,
*\Phi^{\mathfrak{gl}_{\ell+1}}_{\underline{\gamma}} (g)=
\prod_{j=1}^{\ell+1}\,
\pi^{-\frac{\imath\lambda-\imath\gamma_j+\rho_j}{2}}\,
\Gamma\Big(\frac{\imath \lambda-\imath \gamma_{j}+\rho_j}{2}\Big)
\,\,\Phi^{\mathfrak{gl}_{\ell+1}}_{\underline{\gamma}}(g), \ee and
the restriction of $\phi_{\CQ_0(\lambda)}$ to the subspace of
functions satisfying (\ref{equivar}) coincides with  the operator
$\CQ_0^{\mathfrak{gl}_{\ell+1}}(\lambda)$. We shall call such
$\phi_{\CQ_0(\l)}$  a universal Baxter operator.

\begin{te}\label{MainTh}
 Let $\phi_{\CQ_0(\lambda)}(g)$ be a $K$-biinvariant
function on $G=GL(\ell+1,\RR)$  given by \be\label{UBO}
\phi_{\CQ_0(\l)}(g)=2^{\ell+1}|\det g|^{\imath
\lambda+\frac{\ell}{2}} e^{-\pi{\rm Tr} g^tg}. \ee

\noindent i) Then, the action of $\phi_{\CQ_0(\lambda)}$ on the functions
satisfying (\ref{equivar})  descends to the action of
$\CQ_0^{\mathfrak{gl}_{\ell+1}}(\lambda)$ defined  by
(\ref{BAxterredone});

\noindent ii) The action of $\phi_{\CQ_0(\lambda)}$ on modified
Whittaker function
$\Phi^{\mathfrak{gl}_{\ell+1}}_{\underline{\gamma}}(g)$ by a
convolution is given by \be\label{Eigenprop}\bigl(\phi_{\CQ_0(\l)}*
\Phi^{\mathfrak{gl}_{\ell+1}}_{\underline{\gamma}}\bigr)(g)=
L_{\infty}(\lambda)\,\,\Phi^{\mathfrak{gl}_{\ell+1}}_{\underline{\gamma}}(g),
\ee where $L_{\infty}(\lambda)$ is the  local Archimedean $L$-factor
\be L_{\infty}(\lambda)=\prod_{j=1}^{\ell+1}\,
\pi^{-\frac{\imath\lambda -\imath \gamma_j+\rho_j}{2}}
\Gamma\Big(\frac{\imath \lambda-\imath \gamma_{j}+\rho_j}{2}\Big).
\ee
\end{te}
{\it Proof}.

 \noindent {\it i)}. The action of the $K$-biinvariant
function on $\mathfrak{gl}_{\ell+1}$-Whittaker functions is given by
\be\bigl(\phi\,*\,
\Phi^{\mathfrak{gl}_{\ell+1}}_{\underline{\gamma}}\bigr)(g)
=\int_G\, d\tilde{g}\, \phi(g\tilde{g}^{-1})\,
\Phi^{\mathfrak{gl}_{\ell+1}}_{\underline{\gamma}}(\tilde{g})=\int_G\,
 d\tilde{g}\,
\phi(g\tilde{g}^{-1})\,\<k|\pi_{\underline{\gamma}}(\tilde{g})|\psi_R\>.
\ee
Fix the Iwasawa decomposition
$\tilde{g}=\tilde{k}\tilde{a}\tilde{n}$, $\tilde{k}\in K$,
$\tilde{a}\in A$, $\tilde{n}\in N_-$ of a generic element
$\tilde{g}\in G$  and let
$\delta_{B_-}(\tilde{a})=\det_{\mathfrak{n}_-} {\rm
Ad}_{\tilde{a}}$. We shall use  the notation
$d^{\times}a=da\cdot\det(a)^{-1}$ for $a\in A$.   We have for $a\in A$
\bqa\bigl(\phi\,*\,
\Phi^{\mathfrak{gl}_{\ell+1}}_{\underline{\gamma}}\bigr)(a)=
\int_{AN_-}\, d^{\times}\tilde{a}d\tilde{n}\,\delta_{B_-}(\tilde{a})
\,\phi(a\tilde{n}^{-1}\tilde{a}^{-1})\,
\chi_{N_-}(\tilde{n})\,\Phi^{\mathfrak{gl}_{\ell+1}}_{\underline{\gamma}}
(\tilde{a})\,= \nonumber\\
=\int_{A}d^{\times}\tilde{a}\,\, K_{\phi}(a,\tilde{a})\,
\Phi^{\mathfrak{gl}_{\ell+1}}_{\underline{\gamma}} (\tilde{a})\eqa
with \be
K_{\phi}(a,\tilde{a})=\int_{N_-}\!\!d\tilde{n}\,\,\delta_{B_-}(\tilde{a})\,
\,\phi(a\tilde{n}^{-1}\tilde{a}^{-1})\,\chi_{N_-}(\tilde{n}), \\ \nonumber
\chi_{N_-}(\tilde{n})=\exp\Big\{\,2\pi \i
\sum_{i=1}^{\ell}\tilde{n}_{i+1,i}\Big\}.\ee  Thus to prove the first
statement of the Theorem we should prove the following \be
\CQ_0^{\mathfrak{gl}_{\ell+1}}(\underline{x},\underline{y}|\l)
=\int_{N_-}\!\!d\tilde{n}\,\,\delta_{B_-}(\tilde{a})\,
\phi_{\CQ_0(\lambda)}(a\tilde{n}^{-1}\tilde{a}^{-1}|\l)\,\chi_{N_-}(\tilde{n}), \ee
where \be a={\rm\,diag}\,(e^{x_1},\ldots,\,e^{x_{\ell+1}}), \qquad
\tilde{a}={\rm diag}(e^{y_1},\ldots e^{y_{\ell+1}}), \qquad\\
\delta_{B_-}(\tilde{a})=e^{-2\<\rho,\log
\tilde{a}\>}=e^{\sum_{i>j}(y_i-y_j)}. \ee
 For
$g=a\tilde{n}^{-1}\tilde{a}^{-1}$ we have \be
\det\,g=e^{\sum_{i=1}^{\ell+1}(x_i-y_i)}, \qquad
\Tr\,g^tg=\sum_{i=1}^{\ell+1}\,e^{2(x_i-y_i)}+ \sum_{i>j}
u_{ij}^2e^{2(x_i-y_j)}, \ee where $u=\tilde{n}^{-1}\in N_-$.
Taking into account that $\chi_{N_-}(\tilde{n})=\chi_{N_-}(u^{-1})=\exp(-
2\pi \i\sum_{i=1}^{\ell}u_{i+1,i})$ we obtain
\bqa\CQ_0^{\mathfrak{gl}_{\ell+1}}(\underline{x},\underline{y}|\l)=
2^{\ell+1}\int_{N_-}\!\! du\,\,
e^{\sum_{i>j}(y_i-y_j)} \,e^{-2\pi\i\sum_{k=1}^{\ell}u_{i+1,i}}\cdot\\
\nonumber\exp\Big\{\sum_{i=1}^{\ell+1}(\imath\l+\frac{\ell}{2})
(x_i-y_i)- \pi\sum_{i=1}^{\ell+1}e^{2(x_i-y_i)}-
\pi\sum_{i>j}u_{ij}^2e^{2(x_i-y_j)}\Big\}=\\
2^{\ell+1}\exp\Big\{\imath\l\sum_{i=1}^{\ell+1}(x_i-y_i)-
\pi \sum_{i=1}^{\ell+1}e^{2(x_i-y_i)}\Big\}\,e^{\sum_{i>j}(y_i-y_j)}\cdot\\
\nonumber \int_{\mathbb{R}^\ell}\prod_{i=1}^\ell\!du_{i+1,i}\,\,
\exp\Big\{-2\pi \i\sum_{k=1}^{\ell}u_{i+1,i}-\pi\sum_{i=1}^\ell
u_{i+1,i}^2e^{2(x_{i+1}-y_i)}\Big\}\cdot\\ \nonumber
\prod_{i>j+1}\int\!\!du_{ij}\,\, \exp\Big\{-\pi
u_{ij}^2e^{2(x_i-y_j)}\Big\}.\eqa Computing the integrals  by using
the formula\be \int_{-\infty}^{\infty}e^{-\imath\omega
x-px^2}dx=\sqrt{\frac{\pi}{p}}\,\,\,\,e^{\frac{-\omega^2}{4p}}\ee
we readily obtain that \bqa
\CQ_0^{\mathfrak{gl}_{\ell+1}}(\underline{x},\underline{y}|\l)=
\,2^{\ell+1}\exp\Big\{\sum_{i=1}^{\ell+1}(\imath \l+\rho_i)(x_i-y_i)-\\
\nonumber -\pi
\sum_{i=1}^{\ell}\Big(e^{2(x_i-y_i)}+e^{2(y_i-x_{i+1})}\Big)- \pi
e^{2(x_{\ell+1}-y_{\ell+1})}\,\Big\},\eqa where
$\rho_j=\frac{\ell}{2}+1-j,\,\,\,j=1,\ldots,\ell+1$. This completes
the proof of the first statement of the Theorem.

\noindent {\it ii)}. The proof of (\ref{Eigenprop}) follows from the
results of Section 2 $\Box$

It is instructive to provide a direct proof of (\ref{Eigenprop}).
To do so let us first recall  standard facts in  theory of
spherical functions (see \cite{HC} for  details).

There is a general integral expression for the $K$-biinvariant
function in terms of  eigenvalues $\Lambda_{\phi}(\gamma)$
(\ref{phiEigen}). Consider the action on the spherical functions
\be
\varphi_{\gamma}(g)=\<k|\pi_{\gamma}(g)|k\>,
\ee
normalized by
the condition $\varphi_{\gamma}(e)=1$. The explicit integral
representation for $\varphi_{\gamma}(g)$ is
\be
\varphi_{\gamma}(g)=\int_K \,dk\,\, e^{\i\<h(gk),\gamma\>},
\ee
where $\int_K dk=1$ and $h(g)=\log a$ where $g=kan\in KAN_-\to G$ is
the Iwasawa decomposition. Then we have \be
\phi*\varphi_{\gamma}(g)=\Lambda_{\phi}(\gamma)\varphi_{\gamma}(g),
\ee \be \Lambda_{\phi}(\gamma)=\phi*\varphi_{\gamma}(e). \ee Thus
the eigenvalues  can be written  in terms of the spherical transform
as follows   \bqa\label{inteq} \Lambda_{\phi}(\gamma)=\int_G dg
_{}\phi(g^{-1})\varphi_{\gamma}(g)=2^{-(\ell+1)}\int_{A^{+}} \,
d^{\times}a\,\phi(a^{-1})\varphi_{\gamma}(a), \eqa where we have
used the Cartan decomposition $G=KA^+ (M\setminus K)$ to represent
the first integral as an integral over diagonal matrices. Here we
define $A^+=\exp{\frak{a}^+}$, where $\frak{a}^+$ consists of the
diagonal matrices of the form $diag
(e^{x_1},\ldots,e^{x_{\ell+1}}),\,\,\,x_1 < x_2 < \ldots <
x_{\ell+1}$ and $M$ is the normalizer of $\frak{a}$ in $K$. Notice
that $2^{\ell+1}=|M|$.

\begin{prop}
The following integral relation holds \bqa\label{inteqGL}
\Lambda_{\phi_{\CQ_0(\lambda)}}(\gamma)=2^{-(\ell+1)}\int_{A^{+}} \,
d^{\times}a\,\,\phi_{\CQ_0(\lambda)}(a^{-1})\,\varphi_{\gamma}(a)=\\
\nonumber \prod_{j=1}^{\ell+1}\,
\pi^{-\frac{\i\l-\gamma_j+\rho_j}{2}}\,
\Gamma\bigl(\frac{\i\lambda-\imath \gamma_j+\rho_j}{2}\bigr), \eqa
where $\rho_j=\frac{\ell}{2}+1-j,\,\,\,j=1,\ldots,\ell+1$.
\end{prop}
{\it Proof}: Using the integral representation (\ref{inteq}), the
l.h.s. of (\ref{inteqGL})  is given by \be \int_{K\times A^{+}}\, dk
d^{\times}a\,\,\, |\det a|^{-\imath\lambda-\frac{\ell}{2}}e^{-\pi\Tr
(a^ta)^{-1}}\,e^{\i<h(ak),\gamma>}. \ee Using Cartan and Iwasawa
decompositions we have \be\nonumber \int_{K\times
A^{+}}\,dk\,\,d^{\times}a\,\,\, |\det
a|^{-\imath\lambda-\frac{\ell}{2}}e^{-\pi\Tr(a^ta)^{-1}}\,
e^{\i<h(ak),\gamma>}=\ee \be=\int_{K\times A^+\times K} \,dk'\,
d^{\times}a\,dk\,\, |\det
k'ak|^{-\imath\lambda-\frac{\ell}{2}}e^{-\pi {\rm \Tr}
  ((k'ak)^t(kak'))^{-1}}\,e^{\i<h(k'ak),\gamma>}=\ee
\be 2^{\ell+1}\int_{K\times A^+\times M\setminus K} \,dk'\,
d^{\times}a\,dk\,\, |\det
k'ak|^{-\imath\lambda-\frac{\ell}{2}}e^{-\pi {\rm \Tr}
  ((k'ak)^t(kak'))^{-1}}\,e^{\i<h(k'ak),\gamma>}
\ee \be =2^{\ell+1}\int_{G} \,dg\,|\det
g|^{-\imath\lambda-\frac{\ell}{2}}e^{-\pi {\rm \Tr}
(g^tg)^{-1}}\,e^{\i<h(g),\gamma>}=\\ \nonumber
2^{\ell+1}\int_{K\times A\times N_-}\,dn\,d^{\times}a\,dk\,
\delta_{B_-}(a)|\det
a|^{-\i\lambda-\frac{\ell}{2}}e^{-\pi\Tr(n^ta^2n)^{-1}}\,e^{\i\<\log(a),\gamma\>}
\nonumber=\ee \be =2^{\ell+1}\int_{A\times
N_-}\,dn\,d^{\times}a\,\delta_{B_-}(a)|\det
a|^{-\imath\lambda-\frac{\ell}{2}}e^{-\pi\Tr
(n^ta^2n)^{-1}}\,e^{\imath<\log(a),\gamma>}=\\ \nonumber
\prod_{j=1}^{\ell+1}\, \pi^{-\frac{\i\l-\gamma_j+\rho_j}{2}}\,
\Gamma\bigl(\frac{\i\l-\i\gamma_j+\rho_j}{2}\bigr).\ee  Where the
formula
$$\int_{-\infty}^{+\infty}dx e^{\nu
x}e^{-ae^{-2x}}=\frac{1}{2}\,a^{\frac{\nu}{2}}\,\Gamma(-\frac{\nu}{2})$$
was used.
$$\Box$$
The integral operator constructed above can be considered  as a universal
Baxter operator  on matrix elements between the
spherical vector and any other vector in
the representation space. In particular it easy to describe explicitly
an action of the Baxter operators on the space of  zonal spherical
functions. In this case on obtains the  Baxter operator for the  Sutherland
model at a particular value of the coupling constant.

\subsection{Non-Archimedean analog of  Baxter operator}

Let us construct  a non-Archimedean analog of
the universal Baxter $\mathcal{Q}$-operator introduced above.
In the non-Archimedean case the local Hecke algebra
$\mathcal{H}_p=\mathcal{H}(GL(\ell+1,\mathbb{Q}_p),K_p)$,
$K_p=GL(\ell+1,\mathbb{Z}_p)$ is defined
as an algebra of the compactly supported
  $K_p$-biinvariant functions
 on $GL(\ell+1,\mathbb{Q}_p)$. Note that
$K_p$ is a  maximal compact subgroup of $GL(\ell+1,\mathbb{Q}_p)$. Consider a  set
$\{T^{(i)}_p\}$, $i=1,\ldots, (\ell+1)$  of
generators  of $\CH(GL({\ell+1},\mathbb{Q}_p),K_p)$
given by the  characteristic functions of the following  subsets
\bqa
\mathcal{O}_i=K_p \cdot  {\rm diag}(\underbrace{p,\cdots, p}_i
,1\cdots ,1)  \cdot K_p\subset GL({\ell+1},\mathbb{Q}_p).
\eqa
The action of
$T^{(i)}_p$ on functions $f\in C(G/K)$ is then given by the
following integral formula
\bqa
(T^{(i)}_p\,f)\,(g)=\int_{\mathcal{O}_i}\,f(gh)dh.
\eqa
This can be considered as a convolution with characteristic function
$T_p^{(i)}$ of  $\mathcal{O}_i$.
For an appropriately defined non-Archimedean
$\mathfrak{gl}_{\ell+1}$-Whittaker
function $W_{\sigma}$ \cite{Sh}, \cite{CS} one has
\bqa\label{eigenprop1}
T^{(i)}_p\,W_{\sigma}=\,{\rm Tr}_{V_{\omega_i}}\,
\rho_i(\sigma)
\,\,W_{\sigma},
\eqa
where $\rho_i:GL(\ell+1,\CC)\to {\rm End}(V_{\omega_i},\CC)$,
$V_{\omega_i}=\wedge^i \mathbb{C}^{\ell+1}$ is a
representation of $GL(\ell+1,\CC)$ corresponding to the
fundamental  weight $\omega_i$
and $\sigma$ is a conjugacy  class in $GL(\ell+1,\CC)$
corresponding to a non-Archimedean  Whittaker function $W_{\sigma}$.
Note that, in contrast with (\ref{eigenprop1}),
the standard normalization of $T^{(i)}_p$ includes an
additional  factor $p^{-i(i-1)/2}$.
More generally, one considers  Hecke operators $T^{(V)}_p$
associated to arbitrary representations  $\rho_V : GL(\ell+1,\CC)\to
{\rm End}(V,\CC)$ satisfying
\bqa\label{genHeckp}
T^{(V)}_p\,W_{\sigma}=\,{\rm Tr}_{V}\,
\rho_V(\hat{\sigma}_p)
\,\,W_{\sigma}.
\eqa
It is natural to arrange  the generators of $\mathcal{H}_p$ into the
following generating function
\bqa\label{genfuncpad}
T_p(\la)=\sum_{j=1}^{\ell+1} (-1)^j p^{-(\ell+1-j)\la}\, T_p^{(j)}.
\eqa
We introduce  another  generating function
\bqa\label{Baxterpad}
\mathcal{Q}_p^{\mathfrak{gl}_{\ell+1}}(\la)=\sum_{n=0}^{\infty}
\,p^{-n\la}\,\,T^{(S^nV)}_p,
\eqa
where $V=\CC^{\ell+1}$ is the standard representation of $\mathfrak{gl}_{\ell+1}(\CC)$.
The generating functions (\ref{genfuncpad}), (\ref{Baxterpad})
 satisfy  the following relations
\bqa\label{QPADa}
\mathcal{Q}_p^{\mathfrak{gl}_{\ell+1}}(\la)\cdot
\mathcal{Q}_p^{\mathfrak{gl}_{\ell+1}}(\la')
=\mathcal{Q}_p^{\mathfrak{gl}_{\ell+1}}(\la')\cdot
\mathcal{Q}_p^{\mathfrak{gl}_{\ell+1}}(\la),
\eqa
\bqa\label{QPADb}
\mathcal{Q}_p^{\mathfrak{gl}_{\ell+1}}(\la)\cdot T_p(\la')=T_p(\la')\cdot
\mathcal{Q}_p^{\mathfrak{gl}_{\ell+1}}(\la),
\eqa
\bqa\label{QPADc}
1= T_p(\la)\cdot \mathcal{Q}_p^{\mathfrak{gl}_{\ell+1}}(\la),
\eqa
and the operators $T_p(\la)$ and
$\mathcal{Q}_p^{\mathfrak{gl}_{\ell+1}}(\la)$  act on the non-Archimedean
analog of  Whittaker function as
\bqa
T_p(\la)\,\,W_{\sigma}=\det_{V}(1-p^{-\la}\rho_V(\hat{\sigma}_p))\,W_{\sigma},
\eqa
\bqa\label{eigennonar}
\mathcal{Q}^{\mathfrak{gl}_{\ell+1}}_p(\la)\,\,
W_{\sigma}=\det_{V}(1-p^{-\la}\rho_V(\hat{\sigma}_p))^{-1}\,W_{\sigma}.
\eqa
Thus the eigenvalues of
$\mathcal{Q}^{\mathfrak{gl}_{\ell+1}}_p(\la)$
are given by the local non-Archimedean $L$-factors
\bqa\label{eigen1nonar}
L_p(s)=\det_V(1-p^{-s} \rho_V(\hat{\sigma}_p))^{-1},
\eqa
where we use a  more traditional notation $s:=\la$.

Comparing (\ref{QPADa}), (\ref{QPADb}), (\ref{QPADc})
with (\ref{firstpr}), (\ref{secondpr}), (\ref{thirdpr})
one can see that the $\mathfrak{gl}_{\ell+1}$ Baxter $\mathcal{Q}$-operator
appears quite similar to the  generating
function $\mathcal{Q}^{\mathfrak{gl}_{\ell+1}}_p(\la)$ in the Hecke algebra
$\CH(GL(\ell+1,\mathbb{Q}_p),K_p)$ and the  analog of $T_p(\la)$ is
given by (\ref{Genar}). In particular both operators share  the
property that their eigenvalues are given by local $L$-factors.

One can represent  Archimedean and non-Archimedean Baxter operators
in a unified form. Let us rewrite (\ref{Baxterpad}) as
\be\label{BaxnonArch} \CQ^{\mathfrak{gl}_{\ell+1}}(\lambda)(g)
=\sum_{(n_1,,\ldots n_{\ell+1}) \in \mathbb{Z}_+^{\ell+1}}
(p^{n_1}\cdots
p^{n_{\ell+1}})^{\imath\lambda}\delta_{\underline{n}}(g), \ee where
$\underline{n}=(n_1,\ldots, n_{\ell+1})$,
$\delta_{\underline{n}}(g)$ is a characteristic function of
$\CO_{\underline{n}}\subset GL(\ell+1,\mathbb{Q}_p)$ \be
\CO_{\underline{n}}=K_p\cdot {\rm diag} (p^{n_1},\ldots
,p^{n_{\ell+1}})\cdot K_p. \ee On the other hand (universal)
Archimedean Baxter $\CQ$-operator (\ref{UBO}) can be written in the following
form \be \phi_{\CQ_0(\lambda)}(g)=\int dt_1\cdots dt_{\ell+1}
\,(t_1\cdots t_{\ell+1})^{\imath \lambda} e^{-\pi\sum_{j=1}^{\ell+1}
t_i^2} \delta_{\underline{t}}(g), \ee where
$\delta_{\underline{t}}(g)$ is appropriately defined function with the
support at $\CO_{\underline{t}}\subset GL(\ell+1,\mathbb{R})$
\be\label{BaxArch} \CO_{\underline{t}}=K\cdot {\rm diag} (t_1,\ldots
,t_{\ell+1})\cdot K. \ee The integral formulas (\ref{BaxnonArch})
and (\ref{BaxArch}) are compatible in the sense of the standard
correspondence between Archimedean and non-Archimedean integrals
(see e.g. \cite{W}).

\section{Baxter operator for $\mathfrak{so}_{2\ell+1}$}

In the next section we
define a Baxter $\CQ$-operator for
$\mathfrak{g}=\mathfrak{so}_{2\ell+1}$ and
demonstrate that the relation between local $L$-factors and eigenvalues of
$\CQ$-operators holds in this case. More systematic discussion of
the general case will be given elsewhere.

According to  \cite{Ko1},
$\mathfrak{so}_{2\ell+1}$-Whittaker function can be
 written in terms of the invariant
pairing of  Whittaker modules as follows
 \bqa\label{pairingso}
 \Psi^{\mathfrak{so}_{2\ell+1}}_{\lambda} (x)=e^{-
\<\rho,x\>}\<\psi_L\,, \pi_{\lambda}(e^{h_x}) \,\psi_R\>,\, \qquad
\qquad x\in \mathfrak{h}, \eqa  where
$h_x:=\sum\limits_{i=1}^{\ell}\<\omega_i,x\>\,h_i$, $\omega_i$ is a
bases of the fundamental weights of $\mathfrak{so}_{2\ell+1}$. Note
that $\mathfrak{so}_{2\ell+1}$-Whittaker functions are common
eigenfunctions of the complete  set of the commuting
$\mathfrak{so}_{2\ell+1}$-Toda chain Hamiltonians $\CH_{2k}\in{\rm
Diff} (\mathfrak{h})$, $k=1,\cdots, \ell$
 defined by \bqa
\label{hamdef1} \CH^{\mathfrak{so}_{2\ell+1}}_{2k}
\Psi^{\mathfrak{so}_{2\ell+1}}_{\lambda}(x)=e^{-\<\rho,x\>}
\<\psi_L\,,\pi_{\lambda}(e^{h_x})\,c_{2k}\,\psi_R\>,\qquad \eqa
where  $\{c_{2k}\}$ are  generators of the center
$\mathcal{Z}(\mathfrak{so}_{2\ell+1})\subset
\mathcal{U}(\mathfrak{so}_{2\ell+1})$. For the quadratic Hamiltonian
we have \bqa\label{BtwoH}
\CH_2^{\mathfrak{so}_{2\ell+1}}&=&-\frac{1}{2}\sum\limits_{i=1}^{\ell}
\frac{\partial^2}{{\partial x_i}^2}+
\frac{1}{2}e^{x_1}+\sum\limits_{i=1}^{\ell-1} e^{x_{i+1}-x_{i}}
.\eqa Let us introduce a generating function for the
$\mathfrak{so}_{2\ell+1}$-Toda chain Hamiltonians as
\bqa\label{genfunctgl1}
t^{\mathfrak{so}_{2\ell+1}}(\la)=\sum_{j=1}^{\ell}\,(-1)^j
\la^{2\ell+1-2j}\CH_{2j}^{\mathfrak{so}_{2\ell+1}}(x). \eqa Then the
$\mathfrak{so}_{2\ell+1}$-Whittaker function  satisfies  the
following equation \bqa \label{eigeneq1}
t^{\mathfrak{so}_{2\ell+1}}(\la)\,\,\,\Psi^{\mathfrak{so}_{2\ell+1}}_{\underline{\la}}
 (\underline{x})= \la\prod_{j=1}^{\ell}(\la^2-\la_j^2)\,\,\,
\Psi^{\mathfrak{so}_{2\ell+1}}_{\underline{\la}}
(\underline{x}), \eqa
where $\underline{\la}=(\la_1,\cdots ,\la_{\ell+1})$ and
$\underline{x}=(x_1,\ldots,x_{{\ell}+1})$.

\begin{te}
Eigenfunctions of the $\mathfrak{so}_{2\ell+1}$-Toda chain
 admit the integral representation:
 \bqa\label{iterwa}
\Psi^{\mathfrak{so}_{2\ell+1}}
_{\lambda_1,\ldots,\lambda_\ell}(x_{\ell,1},\ldots,x_{\ell,\ell})\,=\,
\int\limits_{\RR^{\ell^2}}
\prod_{k=1}^{\ell-1}\prod_{i=1}^kdx_{k,i}\prod_{k=1}^{\ell}
\prod_{i=1}^kdz_{k,i}\,\, \nonumber e^{{\mathcal
F}^{\mathfrak{so}_{2\ell+1}}(x,z)},\eqa where \bqa \nonumber
\mathcal{F}^{\mathfrak{so}_{2\ell+1}}(x,z)=-\imath\lambda_1(x_{1,1}-2z_{1,1})-
\imath\sum\limits_{n=2}^{\ell}\lambda_n \Big(\sum_{i=1}^nx_{n,i}
-2\sum_{i=1}^{n}z_{n,i}+\sum_{i=1}^{n-1}x_{n-1,i}\Big)-
\\-\Big\{\sum_{n=1}^{\ell}e^{z_{n,1}}+
\sum_{k=2}^{\ell}\sum_{n=k+1}^\ell\Big(e^{x_{n-1,k}-z_{n,k}}+
e^{x_{n,k}-z_{n,k}}\Big)+\\+\sum_{n=k}^\ell\Big
(e^{z_{n,k}-x_{n-1,k-1}}+e^{z_{n,k}-x_{n,k-1}} \Big)
+\sum_{n=1}^{\ell}e^{x_{n,n}-z_{n,n}}\Big\},\nonumber \eqa where we
set $x_i:=x_{\ell,i},\,\,\, 1\le i\leq \ell $.
\end{te}

This integral representation was proposed  in \cite{GLO3} ( we made
an additional  change of variables
$z_{\ell,1}\longmapsto-z_{\ell,1}+
\ln\Big(e^{x_{\ell,1}}+e^{x_{\ell-1,1}}\Big)$ in the integral
representation given in  \cite{GLO3}).

\begin{cor}
The following integral operators
$Q_{\mathfrak{so}_{2\ell-1}}^{\mathfrak{so}_{2\ell+1}}$ provide a
recursive construction of $\mathfrak{so}_{2\ell+1}$-Whittaker
function: \bqa\Psi^{\mathfrak{so}_{2\ell+1}}_{\la_1,\ldots,\la_\ell}
(\underline{x_\ell})=
\int_{\RR^{\ell-1}}\!\!\prod_{i=1}^{\ell-1}d{x_{\ell-1,i}}\,\,
Q_{\mathfrak{so}_{2\ell-1}}^{\mathfrak{so}_{2\ell+1}}
(\underline{x_\ell},\underline{x_{\ell-1}}|\la_\ell)
\Psi^{\mathfrak{so}_{2\ell-1}}_{\la_1,\ldots,\la_{\ell-1}}
(\underline{x_{\ell-1}}),\eqa where \bqa
Q_{\mathfrak{so}_{2\ell-1}}^{\mathfrak{so}_{2\ell+1}}
(\underline{x_\ell},\underline{x_{\ell-1}}|\la_\ell)=
\int_{\RR^{\ell}}\,\,\prod_{i=1}^\ell dz_{\ell,i}\,\,\,
\times  \\
\nonumber\times \exp\Big\{\,-\imath\lambda_\ell\Big(\sum_{i=1}^\ell
x_{\ell,i}-2\sum_{i=1}^{\ell}z_{\ell,i}+
\sum_{i=1}^{\ell-1}x_{\ell-1,i}\Big)\Big\} \,\,\times \\ \nonumber
\times \exp\Big\{-\Big( e^{z_{\ell,1}}+
\sum_{i=1}^{\ell-1}\Big(e^{x_{\ell-1,i}-z_{\ell,i}}+
e^{z_{\ell,i+1}-x_{\ell-1,i}}\Big)+\\ +\nonumber
\sum_{i=1}^{\ell-1}\Big(e^{x_{\ell,i}-z_{\ell,i}}+
e^{z_{\ell,i+1}-x_{\ell,i}}\Big)+e^{x_{\ell,\ell}-z_{\ell,\ell}}\,\Big)
\Big\}.\eqa
For $\ell=1$ we set  \bqa
Q^{\mathfrak{so}_3}_{\mathfrak{so}_1} (x_{1,1};\lambda_1)=\int_{\RR}
dz_{1,1}e^{\imath\lambda_1 x_{1,1}-2\imath\lambda_1
z_{1,1}}\exp\Big\{-\Big(e^{z_{1,1}}+e^{x_{1,1}-z_{1,1}}\Big)\Big\}.\eqa
\end{cor}

Below $\Psi^{\mathfrak{so}_{2\ell+1}}_{\la}(x)$ will
always denote the unique $W$-invariant solution of (\ref{eigeneq1})
(class one principal series Whittaker function). Note that
the space of $W$-invariant Whittaker functions
$\Psi^{\mathfrak{so}_{2\ell+1}}_{\la}(x)$ provides a bases in the
space of $W$-invariant functions on $\mathbb{R}^{\ell}$.

\begin{de} Baxter $\mathcal{Q}$-operator for $\mathfrak{so}_{2\ell+1}$ is given by
\bqa \CQ^{\mathfrak{so}_{2\ell+1}}
 (\underline{y},\underline{x}|\la)=\int_{\RR^{\ell+1}}\,\,\prod_{i=1}^{\ell+1}
dz_{i}\,\exp\Big\{\,-\i\la\Big(\sum_{i=1}^\ell
y_{i}-2\sum_{i=1}^{\ell+1}z_{i}+
\sum_{i=1}^{\ell}x_{i}\Big)\Big\}\times  \\ \nonumber \times
\exp\Big\{-e^{z_{1}}-\sum_{i=1}^{\ell}\Big(e^{y_{i}-z_{i}}+
e^{z_{i+1}-y_{i}}+e^{x_{i}-z_{i}}+
e^{z_{i+1}-x_{i}}\Big)\,\Big\},\eqa where
$\underline{y}=(y_1,\ldots,y_\ell)$ and
$\underline{x}=(x_1,\ldots,x_\ell)$.
\end{de}

\begin{te} Operator $\mathcal{Q}^{\mathfrak{so}_{2\ell+1}}(\la)$
 satisfies the following identities
\bqa\label{firstpr1} \mathcal{Q}^{\mathfrak{so}_{2\ell+1}}(\la)\,
\mathcal{Q}^{\mathfrak{so}_{2\ell+1}}(\la')=
\mathcal{Q}^{\mathfrak{so}_{2\ell+1}}(\la')\,
\mathcal{Q}^{\mathfrak{so}_{\ell+1}}(\la), \eqa
\bqa\label{firstpr15} \mathcal{Q}^{\mathfrak{so}_{\ell+1}}(\la)\cdot
Q^{\mathfrak{so_{2\ell+1}}}_{\mathfrak{so}_{2\ell-1}}(\la')
=\Gamma\Big(\i\la'-\i\la\Big)\Gamma\Big(-\i\la'-\i\la\Big)
Q^{\mathfrak{so_{2\ell+1}}}_{\mathfrak{so}_{2\ell-1}}(\la')\cdot
\mathcal{Q}^{\mathfrak{so}_{\ell-1}}(\la),\eqa \bqa\label{secondpr1}
\mathcal{Q}^{\mathfrak{so}_{2\ell+1}}(\la)\,
T^{\mathfrak{so}_{2\ell+1}}(\la')=
T^{\mathfrak{so}_{2\ell+1}}(\la')\,
\mathcal{Q}^{\mathfrak{so}_{2\ell+1}}(\la), \eqa
\bqa\label{thirdpr1} \la
\mathcal{Q}^{\mathfrak{so}_{\ell+1}}(\la+\imath) =\imath^{2\ell}\,
\mathcal{Q}^{\mathfrak{so}_{\ell+1}}(\la)\, \,
T^{\mathfrak{so}_{2\ell+1}}(\la), \eqa where \bqa\label{Genarso}
T^{\mathfrak{so}_{2\ell+1}}(\underline{x},\underline{y}|\la)=
t^{\mathfrak{so}_{2\ell+1}}(\underline{x},\pr_{\underline{x}}|\la)
\delta^{(\ell)}(\underline{x}-\underline{y}), \eqa \bqa
t^{\mathfrak{so}_{2\ell+1}} (\underline{x},\pr_{\underline{x}}|\la)
=\sum_{j=1}^{\ell+1}(-1)^j
\la^{2\ell+1-2j}\,\CH^{\mathfrak{so}_{2\ell+1}}_{2j}
(\underline{x},\pr_{\underline{x}}). \eqa
\end{te}

{\it Proof}. We will  prove the commutativity of
$\mathcal{Q}$-operators (\ref{firstpr1}). The relation
(\ref{firstpr15}) can be proved using the similar approach. The
other identities  then easily follow.

To prove (\ref{firstpr1}) we should verify the following identity
between the kernels: \bqa\label{commutso}
\int_{\mathbb{R}^{\ell+1}}\,
\mathcal{Q}^{\mathfrak{so}_{2\ell+1}}(\underline{y},\underline{x}|\la)\,\,
\mathcal{Q}^{\mathfrak{so}_{2\ell+1}}(\underline{x},\underline{z}|\la')\,\,
\prod_{j=1}^{\ell+1}d{x_j} =\\\int_{\mathbb{R}^{\ell+1}}\,
\mathcal{Q}^{\mathfrak{so}_{2\ell+1}}(\underline{y},\underline{x}|\la')\,\,
\mathcal{Q}^{\mathfrak{so}_{2\ell+1}}(\underline{x},\underline{z}|\la)\,\,
\prod_{j=1}^{\ell+1}d{x_j}, \eqa where \bqa
\mathcal{Q}^{\mathfrak{so}_{2\ell+1}}
(\underline{y},\underline{x}|\la)=\int_{\RR^{\ell+1}}\,\,\prod_{i=1}^{\ell+1}
du_{i}\,\exp\Big\{\,-\i\la\Big(\sum_{i=1}^\ell
y_{i}-2\sum_{i=1}^{\ell+1}u_{i}+
\sum_{i=1}^{\ell}x_{i}\Big)\Big\}\times  \\ \nonumber \times
\exp\Big\{-e^{u_{1}}-\sum_{i=1}^{\ell}\Big(e^{y_{i}-u_{i}}+
e^{u_{i+1}-y_{i}}+e^{x_{i}-u_{i}}+
e^{u_{i+1}-x_{i}}\Big)\,\Big\},\eqa \bqa
 \mathcal{Q}^{\mathfrak{so}_{2\ell+1}}
 (\underline{x},\underline{z}|\la')=
\int_{\RR_{\ell+1}}\,\,\prod_{i=1}^{\ell+1}
dv_{i}\,\exp\Big\{\,-\i\la'\Big(\sum_{i=1}^\ell
x_{i}-2\sum_{i=1}^{\ell+1}v_{i}+
\sum_{i=1}^{\ell}z_{i}\Big)\Big\}\times  \\ \nonumber \times
\exp\Big\{-e^{v_{1}}-\sum_{i=1}^{\ell}\Big(e^{x_{i}-v_{i}}+
e^{v_{i+1}-x_{i}}+e^{z_{i}-v_{i}}+
e^{v_{i+1}-z_{i}}\Big)\,\Big\}.\eqa The proof is given by the
following  sequence of  elementary transformations. Let us first
make a change of variables $u_i$ and $v_i$ in (\ref{commutso}): \be
u_1\longmapsto-u_1+\ln\Big(e^{y_1}+e^{x_1}\Big), &\\ \nonumber
u_i\longmapsto-u_i-\ln\Big(e^{y_{i-1}}+e^{x_{i-1}}\Big)+
\ln\Big(e^{y_i}+e^{x_i}\Big), & 1<i\leq\ell,\ee
\be v_1\longmapsto-v_1+\ln\Big(e^{x_1}+e^{z_1}\Big), &\\
\nonumber v_i\longmapsto-v_i-\ln\Big(e^{x_{i-1}}+e^{z_{i-1}}\Big)+
\ln\Big(e^{x_i}+e^{z_i}\Big), & 1<i\leq\ell.\ee
We introduce additional integration variables $u_{\ell+1}$ and
$v_{\ell+1}$ in (\ref{commutso})
using integral formulas:
 \be\Big(e^{-y_\ell}+e^{-x_\ell}\Big)^{-2\i\la'}=\Gamma(2\i\la')^{-1}
\int_{\RR}\,du_{\ell+1}\,\exp\Big\{2\i\la'u_{\ell+1}-
e^{u_{\ell+1}-y_\ell}-e^{u_{\ell+1}-x_\ell}\Big\},\\
\Big(e^{-x_\ell}+e^{-z_\ell}\Big)^{-2\i\la}= \Gamma(2\i\la)^{-1}
\int_{\RR}\,dv_{\ell+1}\,\exp\Big\{2\i\la'v_{\ell+1}-
e^{v_{\ell+1}-x_\ell}-e^{v_{\ell+1}-z_\ell}\Big\}.\ee
Then let us modify the variables $x_i,\,i=1,\ldots,\ell$ as: \bqa
x_i\longmapsto-x_i-\ln\Big(e^{-u_i}+e^{-z_i}\Big)+
\ln\Big(e^{u_{i+1}}+e^{z_{i+1}}\Big).\eqa
and use the following integral representations to introduce additional
variables $x_0$ and $x_{\ell+1}$
\bqa \Big(e^{u_1}+e^{v_1}\Big)^{-\i(\la+\la')}&=&\Gamma(\i(\la+\la'))^{-1}
\int_{\RR}\, dx_0\,\,\exp\Big\{-\i(\la+\la')x_0-
e^{u_1-x_0}-e^{z_1-x_0}\Big\},\nonumber \\
\Big(e^{-u_{\ell+1}}+e^{-v_{\ell+1}}\Big)^{\i(\la+\la')}&=&
\Gamma(-\i(\la+\la'))^{-1}\times
\nonumber \\\nonumber &&\times \int_{\RR}
\,dx_{\ell+1}\,\exp\Big\{-\i(\la+\la')x_{\ell+1}-
e^{x_{\ell+1}-u_{\ell+1}}-e^{x_{\ell+1}-v_{\ell+1}}\Big\}.\eqa
Now we make the following sequence of
changes of the variables: \be
u_1\longmapsto-u_1-\ln\Big(1+e^{-x_0}\Big)+
\ln\Big(e^{y_1}+e^{x_1}\Big), &\\
\nonumber u_i\longmapsto-u_i-\ln\Big(e^{y_{i-1}}+e^{x_{i-1}}\Big)+
\ln\Big(e^{y_i}+e^{x_i}\Big), & 1<i\leq\ell,\\
u_{\ell+1}\longmapsto-u_{\ell+1}+x_{\ell+1}-
\ln\Big(e^{-y_\ell}+e^{-x_\ell}\Big), &\ee

\be
v_1\longmapsto-v_1-\ln\Big(1+e^{-x_0}\Big)+
\ln\Big(e^{x_1}+e^{z_1}\Big), &\\
\nonumber v_i\longmapsto-v_i-\ln\Big(e^{x_{i-1}}+e^{z_{i-1}}\Big)+
\ln\Big(e^{x_i}+e^{z_i}\Big), & 1<i\leq\ell,\\
v_{\ell+1}\longmapsto-v_{\ell+1}+x_{\ell+1}-
\ln\Big(e^{-x_\ell}+e^{-z_\ell}\Big), &\ee

\be
x_0\longmapsto-x_0+\ln\Big(e^{u_1}+e^{z_1}\Big), &\\
x_i\longmapsto-x_i-\ln\Big(e^{-u_i}+e^{-z_i}\Big)+
\ln\Big(e^{u_{i+1}}+e^{z_{i+1}}\Big), & 1\leq i\leq\ell,\\
x_{\ell+1}\longmapsto-x_{\ell+1}-
\ln\Big(e^{-u_{\ell+1}}+e^{-z_{\ell+1}}\Big).&\ee

One integrates out the variables $x_0$ and $x_{\ell+1}$
and modifies the  variables $u_i$ and $v_i$ as follows  \be
u_1\longmapsto-u_1+\ln\Big(e^{y_1}+e^{x_1}\Big), &\\
\nonumber u_i\longmapsto-u_i-\ln\Big(e^{y_{i-1}}+e^{x_{i-1}}\Big)+
\ln\Big(e^{y_i}+e^{x_i}\Big), & 1<i<\ell,\\
u_{\ell}\longmapsto-u_{\ell}-
\ln\Big(e^{-y_{\ell-1}}+e^{-x_{\ell-1}}\Big), &\ee

 \be
v_1\longmapsto-v_1+\ln\Big(e^{x_1}+e^{z_1}\Big), &\\
\nonumber v_i\longmapsto-v_i-\ln\Big(e^{x_{i-1}}+e^{z_{i-1}}\Big)+
\ln\Big(e^{x_i}+e^{z_i}\Big), & 1<i<\ell,\\
v_{\ell}\longmapsto-v_{\ell}-
\ln\Big(e^{-x_{\ell-1}}+e^{-z_{\ell-1}}\Big). &\ee
Integrating out  $u_{\ell+1}$ and $v_{\ell+1}$,
one completes the proof of (\ref{firstpr1}) $\Box$

\begin{cor} The following identity  holds
\bqa\label{BIntertw}
\int_{\RR^{\ell}}\,\prod_{i=1}^{\ell}dx_{\ell,i}
\CQ^{\mathfrak{so}_{2\ell+1}}(\underline{y},\underline{x}|\gamma)\,
\Psi^{\mathfrak{so}_{2\ell+1}}_{\underline{\la}}(\underline{x})=
\prod_{i=1}^\ell\Gamma\Big(\i \la_i-\i\gamma \Big)
\prod_{i=1}^\ell\Gamma\Big(-\i \la_i-\i \gamma \Big)\,\,
\Psi^{\mathfrak{so}_{2\ell+1}}_{\underline{\la}}(\underline{y}).\eqa
\end{cor}
Finally let us note  that this result is in agreement with the interpretation of the
eigenvalues of $\mathcal{Q}$-operators as local Archimedean
$L$-functions corresponding to automorphic representations of
reductive Lie groups discussed above.

\end{document}